\documentclass[11pt,a4paper,twoside]{amsart}
\usepackage[latin1]{inputenc}
\usepackage[T1]{fontenc}
\usepackage{eucal,latexsym,amssymb,amscd,a4wide,pifont,pstricks,pst-node}

\theoremstyle{plain}
\newtheorem{thm}{\bfseries Theorem}[section]

\newtheorem{lem}[thm]{\bfseries Lemma}

\newtheorem{sublem}[thm]{\bfseries Sub-Lemma}
\newtheorem{prop}[thm]{\bfseries Proposition}
\newtheorem{crit}[thm]{\bfseries Criterion}

\newtheorem{cor}[thm]{\bfseries Corollary}

\newtheorem{df}[thm]{\bfseries Definition}
\theoremstyle{remark}
\newtheorem{ex}[thm]{\bfseries Example}
\newtheorem{nota}[thm]{\bfseries Notation}

\newtheorem{obs}[thm]{\bfseries Observation}
\newtheorem{rem}[thm]{\bfseries Remark}

\numberwithin{equation}{section}

\newenvironment{enum}{\begin{dingautolist}{202}}{\end{dingautolist}}
\DeclareMathSymbol{\Z}{\mathalpha}{AMSb}{"5A} 
\DeclareMathSymbol{\PP}{\mathalpha}{AMSb}{"50} 
\DeclareMathSymbol{\Q}{\mathalpha}{AMSb}{"51}
\DeclareMathSymbol{\N}{\mathalpha}{AMSb}{"4E}
\DeclareMathSymbol{\R}{\mathalpha}{AMSb}{"52}
\newcommand{\A}{\mathcal{A}} 
\newcommand{\RR}{\mathcal{R}} 

\newcommand{\wt} {\widetilde}
\newcommand{\crr} {{\rm cr}}

\newcommand{\F}{\mathbb{F}}

\newcommand{\fl}{\longrightarrow}
\newcommand{\iso}{\cong}

\newcommand{\U}[1]{#1 ^{\times}} 

\newcommand{\pol} {\mathcal D}
\newcommand{\C} {\mathbb{C}}

\newcommand{\CC} {\mathcal{C}}
\newcommand{\B} {\mathcal B}

\newcommand{\CR} {\mathcal R}
\newcommand{\sterling} {£}

\newcommand{\Wedge}{\bigwedge}

\newcommand{\doubledot}{{\scriptscriptstyle\bullet\bullet}}
\newcommand{\floor}[1]{\left\lfloor #1 \right\rfloor}

\newcommand{\ourbeta}{\mathfrak{b}} 

\newcommand{\scg}[1]{\mathfrak{p}(#1)} 

\newcommand{\newsubsubsection}[1]{\subsubsection{{\bfseries\itshape #1}}}
\begin{document}
\title{On Poly(ana)logs I}
\author{Philippe Elbaz-Vincent}
\thanks{The first author was partly supported by a Marie Curie fellowship of the EU}
\author{Herbert Gangl}
\thanks{The second author is supported by a Habilitationsstipendium der Deutschen Forschungsgemeinschaft} 
\begin{abstract} 
We investigate a connection between the differential of polylogarithms
(as considered by Cathelineau) and a finite variant of them. This allows
to answer a question raised by Kontsevich concerning the construction of
functional equations for the finite analogs, using in part the $p$-adic version of polylogarithms and recent work of Besser.   
Kontsevich's original unpublished note is supplied (with his kind
permission) in an "Appendix" at the end of the paper.
\end{abstract}
\maketitle
\tableofcontents

\part*{Introduction and Motivation}
In an unpublished note \cite{kontse} (included as an Appendix) Kontsevich defined the ``$1\frac12$-logarithm'', associated to a prime $p$,  as the
truncated power series of $-\log(1-x)$ (for which we propose the ``truncated''
letter $\sterling$, pronounced ``sterling'') as a function from $\Z/p$ to $\Z/p$:
   $$ \sterling_1(x)=\sterling^{(p)}_1(x)=\sum_{k=1}^{p-1} \frac{x^k}{k}\quad
   \pmod p.$$

For reasons which become apparent below we refer to it as the {\sl finite 1-logarithm}.
Kontsevich observed that it satisfies a functional equation which is known in the
literature as the fundamental equation of information theory (see \cite{ad}), and
provided a cohomological interpretation of the equation.

Cathelineau \cite{cathe-pol} was led to the same equation by considering an ``infinitesimal'' version of a 
one-valued cousin of the dilogarithm function which is defined over $\C$. He had encountered the fundamental equation of information
theory already in \cite{cathe-sl2} where, motivated by questions arising from Hilbert's third problem, he
deduced an infinitesimal version of the famous Bloch--Suslin complex (which calculates certain algebraic $K$-groups
of a field). Furthermore, he provided a homological interpretation of the equation. Cathelineau extended his
results to infinitesimal versions of higher polylogarithms, and in particular---by mimicking Goncharov's setup 
\cite{g1} which generalizes the Bloch--Suslin complex---deduced an infinitesimal analogue of Goncharov's complexes.
In the process, he produced the generic functional equation for the infinitesimal trilogarithm which contains
22 terms in 3 variables.

Kontsevich had asked explicitly in \cite{kontse} for functional equations similar to the fundamental equation of 
information theory for the next case, i.e.~for the case of the finite dilogarithm $\sterling_2(x)=\sum_{k=1}^{p-1} x^k/k^2$.  Guided by the analogy between finite 1-logarithm and the infinitesimal dilogarithm, it was found that 
Cathelineau's equation for the infinitesimal trilogarithm is also satisfied by $\sterling_2$ and
provides an answer to Kontsevich's question. Furthermore, $\sterling_2$ is characterized by the latter
equation (actually, it is already characterized by certain specializations).

In fact we get a stronger statement: {\it each of the functional equations for the infinitesimal $n$-logarithm in this paper---and this includes the
distribution formulas for any $n$---has been proved
for the {\sl finite $(n-1)$-logarithm} (whose definition should be clear by the above)}.

What is more, there is a whole machinery to obtain this type of functional equations: on the one hand,
Cathelineau had given a tangential procedure for elements in $\Z[F]$ (for certain fields $F$)
which is compatible with the passage from functional equations for the dilogarithm to equations for 
the infinitesimal dilogarithm. It turns out (see §\ref{deriving}) that the same is true for higher polylogarithms, and we will show how we can get a functional equation for an  infinitesimal $n$-logarithm by ``taking the derivative'' of a functional equation for the classical $n$-logarithm relatively to an absolute derivation over $F$.
On the other hand, since $p$-adic polylogarithms in the sense of Coleman \cite{coleman} satisfy the same 
functional equations as the classical ones by work of Wojtkowiak \cite{wojt2} (for a more precise statement
cf.~\S\ref{red_mod_p}), one arrives via Cathelineau's tangential procedure (proved by him in characteristic 0)
at its $p$-adic equivalent and one could hope that there is a version of $p$-adic polylogarithms whose 
appropriate differential reduces to the finite polylogarithms.
This hope (vaguely anticipated in  \cite{pev-essen}) has been made precise by Kontsevich (private communication) and was subsequently proved (in a slightly 
modified form) by Besser \cite{Besser}. 
Combining the above, we obtain a recipe for deducing functional equations for $\sterling_{n-1}$
from functional equations for the $n$-logarithm, and thus we get analogues of distribution relations for 
each $n$ and further ``non-trivial'' ones at least up to $n=7$ (cf.~\cite{zag-app}, \cite{gangl1}).
The properties stated motivate the terminology of ``poly(ana)logs'' for the different analogues of polylogs. 
To help the reader to understand the interdependencies between the notions already discussed, we give the following picture, which can serve as a guideline for the paper:

\newpage
\centerline{\small\bf The conceptual relationship between the different Poly(ana)logs}\quad\\[2pt]
\psset{arrows=->,labelsep=5pt,linecolor=gray,framearc=.2}

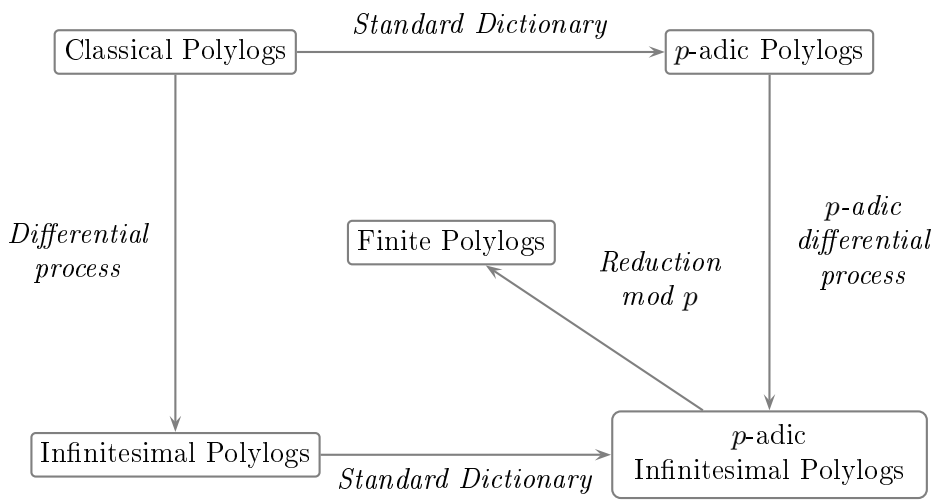
\begin{center}
$\begin{array}{ccc}
\rnode{a}{\psframebox{\textrm{Classical Polylogs}}} &  & \quad\rnode{b}{\psframebox{\textrm{$p$-adic Polylogs}}}\\[2cm]
& \rnode{c}{\psframebox{\textrm{Finite Polylogs}}} & \\[2cm]
\rnode{d}{\psframebox{\textrm{Infinitesimal Polylogs}}} &  & \quad\rnode{e}{\psframebox{\begin{tabular}{c}\textrm{$p$-adic} \\ \textrm{Infinitesimal Polylogs}\end{tabular}}}
\ncline{a}{b}\Aput{\textrm{\emph{Standard Dictionary}}}
\ncline{a}{d}\Bput{\textrm{\emph{\begin{tabular}{c}Differential \\ process\end{tabular}}}}
\ncline{b}{e}\Aput{\textrm{\emph{\begin{tabular}{c}$p$-adic\\ differential \\ process\end{tabular}}}}
\ncline{e}{c}\Bput{\textrm{\emph{\begin{tabular}{c}Reduction \\ mod $p$\end{tabular}}}}
\ncline{d}{e}\Bput{\textrm{\emph{Standard Dictionary}}}
\end{array}$\\[1.3cm]
\end{center}

The present paper investigates the basic properties of the infinitesimal version of polylogarithms, including the $p$-adic ones, and their relationship with the finite polylogarithms and also with the classical polylogarithms via the ``derivation map'' (section \ref{deriving}). In particular, the answer to Kontsevich's question can be found in section 4 (Theorem \ref{L2.1}), together with a proof of the unicity of $£_2$ (Theorem \ref{unicity_of_sterling2}). The sequel paper \cite{polyana2} exhibits 
interrelationships among the polylogarithmic groups and also among their infinitesimal versions, introduces finite
versions of the so-called ``multiple polylogarithms'' (cf.~e.g.~\cite{gonmulpol}) 
and in particular some multiplicative structure related to them:
it turns out that the proofs of the identities for the finite field case are far from trivial, and
especially the most conceptual one found for Cathelineau's 22-term equation involves an identity expressing
$\sterling_1(a)\sterling_1(b)$ in terms of $\sterling_2$ only.
The special case of $a=b$ in the latter product is an identity found by Mirimanoff which is crucial for proving his
criteria for Fermat's last theorem---the finite polylogarithms have appeared in the literature prominently 
in the guise of ``Mirimanoff polynomials'' (cf.~Ribenboim's 13 Lectures \cite{rib}).
Others of Mirimanoff's identities can be reinterpreted in terms of functional equations of finite 
polylogarithms (actually, ``multiple polylogarithms'') which might nurture the hope that further knowledge 
concerning the latter could provide more obstacles for a solution of FLT to exist (but this may well turn out to be
a too pollyanna\footnote{{\bf Pollyanna.} The name of the heroine of stories written by 
Eleanor Hodgman Porter (1868-1920), American children's author, used with allusion to her skill at
the `glad game' of finding cause for happiness in the most disastrous
situations; {\bf one who is unduly optimistic or achieves happiness through
self-delusion}. [Oxford English Dictionary~2]} attitude)...\\

The organisation of the present work is as follows:\\[2pt]
Part I is dedicated to the introduction of classical and infinitesimal polylogarithms (in characteristic 0) and their associated functional equations and groups. In particular we re-introduce several notions of Cathelineau \cite{cathe-sl2,cathe-pol} and give  complementary properties.\\
Part II introduces the finite polylogs, the functional equations that they satisfy and give their characterizations (section \ref{finite}). We also introduce in the section \ref{deriving} the construction of the ``derivation map'' and show that functional equations for classical polylogs give rise to functional equations for infinitesimal polylogs. The last section of this part (section \ref{red_mod_p}) introduces the $p$-adic methods, and shows (Corollary \ref{surprise}), via Besser's result, that functional equations for infinitesimal $p$-adic polylogs produce functional equations for finite polylogs (under mild assumptions).\\
Finally, the main proofs of Part II are given in Part III.\\
The paper ends with a reproduction of the note of Kontsevich \cite{kontse}, originally written for a private booklet dedicated to Friedrich Hirzebruch on the occasion of his ``Emeritierung'' (retirement). We are grateful to him for letting us include it
as an appendix.\\[3pt]

\noindent{\Small {\bf Acknowledgements} : We would like to express our sincere 
gratitude for financial support
to DFG, EU (Marie Curie fellowship program), and to the   following 
institutions for their hospitality: the Laboratoire Dieudonn\'e UMR CNRS 6601 of the University of Nice--Sophia Antipolis where the whole project started out (and where the second author made several visits working on this project),  the  Institut für Experimentelle Mathematik (Essen), the FB6 Mathematik of the Universit\"at-GH Essen, the I.H.\'E.S. and the MPI f\"ur Mathematik Bonn. We want to give our ``Herzliches Dankesch\"on'' to G. Faltings and G. Harder for their invitation to the MPI Conference on Polylogarithms at Schloss Ringberg. The first author wants to thank the Mathematical Institute of the University of Lausanne (and especially D. Arlettaz) where he was working on the subject. The second author wants to thank the Laboratoire G.T.A. UMR CNRS 5030 of the University Montpellier II for its hospitality while he was visiting.
During the long period of gestation of this work we had the opportunity to discuss it with several people. It is a pleasure to thank them here: we are grateful to H. Esnault, G. Frey, A. Goncharov, G. Mersmann, J. Nekova\v r and Z. Wojtkowiak (and also to D. Grayson for pointing out the origin of
Pollyanna). We particularly want to thank A. Besser for sending us a preliminary version of his note and for helpful remarks, P. Colmez for his patience explaining to us several $p$-adic features, C. Soul\'e for his constant interest and questions in our work, as well as D. Zagier for enlightening and stimulating comments. Last but not least, we want to express our warmest thanks to the (inadvertent) initiators of the story: J.-L. Cathelineau and M. Kontsevich, for their sustained encouragement and propulsive discussions, and  without whom this work would not exist.}



\part{Preliminary Background}

\section{Definitions of polylogarithms and their analogues (in characteristic 0)}
In the following we will recall some standard, and some less standard, facts about
polylogarithms and their functional
equations. The main references will be Zagier \cite{z1} and Goncharov
\cite{g3} (for the classical case) as well as Cathelineau \cite{cathe-pol} (for the 
infinitesimal case).
\subsection{Classical and one-valued Polylogarithms}\label{classical-poly}
Let $\,n\geqslant 1\,$, and   $\,{\pol }_{n} : \C\rightarrow{\R}(n-1)\,$ be the
 Bloch/Wigner/Ramakrishnan/Zagier/Wojtkowiak function \cite{z1,g3,cathe-pol,wojt1}, or {\sl modified $n$th polylogarithm}, defined by 
$${\pol }_{n}(z)=\Re_{n} 
\left(\sum_{k=0}^{n-1}\frac{2^{k}B_{k}}{k!}{\rm log}^{k}|z|Li_{n-k}(z)\right),$$
where $\,{\Re}_{n}\,$ denotes ${\rm Re}$ or $i{\rm Im}$, and  
$\,{\R}(n)={\R}\,$ or $\,i{\R}\,$, 
depending on whether $\,n\,$ is even or odd. The $\,B_{k}\,$ are the Bernoulli numbers ($\,B_0=1\,$,
$\,B_1=-\frac12\,$, $\,B_2=\frac16\,$, $\,B_3=0\,$, \dots),
and  $\,Li_{m}\,$ denotes the classical  $\,m$-logarithm
$$ Li_m(z)=\sum_{n=1}^\infty \frac{z^n}{n^m}\,,\qquad |z|<1,$$
which can be analytically continued to the cut plane $\,\C-[1,\infty)\,$ 
\cite{z1}.
For example, we have,
$${\pol }_{1}(z)=-\log|1-z|,$$
$${\pol }_{2}(z)=i\,{\rm Im}\bigg(Li_{2}(z) + \log(1-z)\log |z|\bigg),$$
$${\pol }_{3}(z)={\rm Re}\Big(Li_{3}(z) - \log|z|Li_{2}(z) - \frac{1}{3}
\log^{2}|z|\log(1-z)\Big).$$

\begin{rem} \begin{enumerate}
\item The virtue of these modifications of classical polylogarithms lies
in the fact that they are one-valued functions on the whole complex plane (at the points $0$ and $1$ they are defined by continuity)---as opposed to the multi-valued 
classical polylogarithm functions---and that they satisfy ``clean'' functional
equations (i.e. without lower order terms such as products of polylogarithms of lower degrees).\\
\item Instead of the above $\pol_n$ there is also the closely related real-valued function $P_n$ (originally introduced 
by Zagier \cite{z1}) widely used, and also denoted ${\mathcal L}_n$, e.g. \cite{g1}. It differs from
$\pol_n$ only by a possible factor of $i$.\\
\item {\it Polylogarithms of a real variable}. In a similar manner one can define real valued functions  
as given by Zagier \cite{z1} (eq. (31), p.412), cf. also Lewin \cite{lewin3}
(eq.(16), p.7), which could be called {\it Rogers polylogarithms} in view of Rogers's 
investigations in the case $\,n=2\,$ \cite{Rog}:
$$ L_n(x)=\sum_{j=0}^{n-1}\frac{(-\log|x|)^j}{j!}\,Li_{n-j}(x)\ +\ 
\frac{(-\log|x|)^{n-1}}{n!}\,\log|1-x|\,,\quad |x|\leqslant 1,$$
and for $\,|x|>1\,$ via the inversion relation
$$L_n\Big(\frac1x\Big)=(-1)^{n-1}L_n(x)\,.$$
\end{enumerate}
\end{rem}

\subsection{Infinitesimal polylogarithms}
We mainly follow the presentation in Cathelineau \cite{cathe-pol}. Differentiating the
functions $\,{\pol }_{n}\,$ gives (see \cite{cathe-pol}, p.\,1328)
$$\frac{\partial}{\partial z}
{\pol }_{n}(z) = -\sum_{k=1}
^{n-1}\frac{2^{k-1}B_{k}}{k!}\frac{{\rm log}^{k-1}|z|}{z}{\pol }_{n-k}(z) 
+\frac{2^{n-2}B_{n-1}}{(n-1)!}\frac{{\rm log}^{n-1}|z|}{1-z},$$
and

$$ \frac{\partial}{\partial \bar z}
{\pol }_{n}(z) = (-1)^{n-1}\overline{\frac{\partial}{\partial z}
{\pol }_{n}(z)}.$$
Finally we can deduce the expression for $\,d{\pol }_{n}(z)\,$
\begin{align*}
d{\pol }_{n}(z) &= \frac{\partial}{\partial z}
{\pol }_{n}(z)dz + \frac{\partial}{\partial \bar z}
{\pol }_{n}(z)d\bar z \\
 &= -\sum_{k=1}^{n-1}\left(
\frac{2^{k}B_{k}}{k!}\,{\rm log}^{k-1}|z|\,{\pol }_{n-k}(z)\Upsilon_{k}(z)
\right)
-\frac{2^{n-1}B_{n-1}}{(n-1)!}\,{\rm log}^{n-1}|z|\,\Upsilon_{n-1}(1-z).
\end{align*}
If $\,k\,$ is even : $\,\Upsilon_{k}(z)=d\,{\rm log}|z|\,$, and if $\,k\,$ is odd :  
$\,\Upsilon_{k}(z)=di{\rm arg}(z)\,$. 
The main examples are

$$d{\pol }_{1}(z) = - d\log|1-z|,$$
$$d{\pol }_{2}(z) = -\log|1-z|\,di{\rm arg}(z) +\log|z|\,di{\rm arg}(1-z),$$
$$d{\pol }_{3}(z) = {\pol }_{2}(z)di{\rm arg}(z) + \frac{1}{3}\log|z|\,
\big(\log|1-z|\,d\log|z|- \log|z|\,d\log|1-z|\big).$$

\begin{rem}Goncharov \cite{g1}(Prop.~1.18) had deduced a slightly 
different, but equivalent, 
formula earlier (the terms which seem a priori different---he wrote $\,d\,{\rm log}|z|\,$
instead of $\,d\,{\rm arg}(z)\,$---turn out to be multiplied by a Bernoulli 
number $\,B_k\,$ which is zero since $\,k\,$ is odd).
\end{rem}

\section{Groups related to polylogarithms}\label{polygroups}
\def \Fdotdot {F^{\scriptscriptstyle\bullet\bullet}}
In the following, $F$ will denote a field, and we abbreviate  $\Fdotdot=F-\{0,1\}$. We can think of it as 
a doubly punctured affine line over $F$.
\subsection{The scissors congruence group} We define the scissors congruence group $\scg{F}$ as the quotient of $\Z[\Fdotdot]$ by the subgroup generated by the elements 
$$\left[ a \right] \,-\, \left[ b \right] \,+\,
\left[ \frac{b}{a} \right] \,-\,\left[
\frac{1-b}{1-a} \right] \,+\,\left[
\frac{1-b^{-1}}{1-a^{-1}} \right], $$
whenever such an expression makes sense. The relation is the  famous five term equation for the dilogarithm
(first stated by Abel, cf.~\cite{lewin1}). This group, which has a geometric origin (see for instance \cite{dupont/sah}), captures the algebraic properties of the dilogarithm, more precisely one has
 \begin{prop} If $F\subset \C$, then the dilogarithm $\pol_2$  is defined on $\scg{F}$. 
\end{prop}
Suslin's definition of the Bloch group of a field is given by the following exact sequence (see \cite{suslin}),
\begin{equation}0 \to B(F) \to \scg{F} \overset{\lambda}{\to} (\U{F} \otimes_\Z \U{F})_s \to K^M_2(F) \to 0,\label{B2}\end{equation}
where $K^M_2(F)$ is the Milnor $K_2$ of the field $F$ (see \cite{rosenberg}, chapter 4), $(\U{F} \otimes_\Z \U{F})_s$ is the quotient of $\U{F} \otimes_\Z \U{F}$ by the subgroup generated by the elements of the kind $x \otimes y+y \otimes x$.  The map $\lambda$ is then defined by $\lambda([a])=a\otimes (1-a)$ and the {\it Bloch group} of $F$ is defined as the kernel of this map. 
\begin{rem}\begin{enumerate}
\item In \cite{dupont/sah}, Dupont and Sah have studied in detail the scissors congruence group and also its connection to the dilogarithm.\\
\item If $F$ is an infinite field, the precise relationship between $K_3(F)$ and $B(F)$ is described by Suslin in \cite{suslin}, and rationally we have  $K_3(F)_\Q\iso B(F)_\Q$ 
which gives a description of $K_3(F)_\Q$ in terms of generators and relations.\\

\item Weibel \cite{weibel} has computed the group $B(F)$ if $F$ is a finite field and has shown that it has the same relationship to $K_3$ as in the case of infinite fields.\\
\item The original definition of Bloch \cite{bloch1}(Lecture 6, p.59) is given by the following exact sequence
$$0 \to \B (F)\to  \A (F) \overset{\lambda}{\to}\U{F} \otimes_\Z \U{F}  \to K_2(F) \to 0,$$
where $\A(F) $ is just the group $\Z[\Fdotdot]$. Notice that he also generalized the definition to rings in order to prove some rigidity property \cite{bloch1}(pp.62-68). Moreover, he obtained a map between $\B(F)$ and $K_3^{ind}(F)/Tor_1^\Z(\U{F},\U{F})$ for any algebraically closed field $F$ \cite{bloch1}(pp.71-72). (Here, 
$K_3^{ind}(F)$ denotes the quotient of $K_3(F)$ by the image of $K_3^M(F)$
in $K_3(F)$.) Later, Suslin \cite{suslin} showed that we have an analogous map with $B(F)$ and that, modulo 2-torsion, this map is an isomorphism.\\
\item In fact the exact sequence (\ref{B2}) holds also for ``rings with many units'', such as semilocal rings with infinite residue fields (this is a consequence of results in \cite{pev}).
\end{enumerate}
\end{rem}

\subsection{Polylogarithmic groups and Goncharov complexes}\label{goncomplex} Zagier has generalized in \cite{z1} (section 8) the construction of the Bloch group to higher $n$ and defined
{\sl higher Bloch groups}, on which the corresponding polylogarithm functions $\,\pol_n\,$ are defined.
They are constructed by an 
inductive procedure which has been made more conceptual by Goncharov 
whose framework we adopt here. Let $\PP^1(F)$ be the projective line over $F$. The construction of an
intermediate group $\B_n(F)$, descriptively called {\sl polylogarithmic group} in \cite{jlc}, proceeds by induction on $n\geqslant 2$. We first need to construct certain subgroups $\A_n(F)$ and $\RR_n(F)$ of $\Z[\PP^1(F)]$. Suppose that $\RR_n(F)$ is defined, then we set 
$$\B_n(F)=\Z[\PP^1(F)]/\RR_n(F).$$

Define the morphisms

$$\begin{array}{cc}
\delta_{2} =\delta_{2,F} : & {\Z}[{\PP}^{1}(F)]\rightarrow 
\dfrac{\bigwedge^{2}_\Z F^{\times}}{  (2-torsion)}, \\
 & [x]\mapsto \left\{
\begin{tabular}{ll}
0 & if $x= 0,1,\infty$,\vphantom{\Big|} \\
$(1-x)\wedge x$ & otherwise,
\end{tabular}
\right.
\end{array}$$
and for $n\geqslant 3$
$$\begin{array}{cc}
\delta_{n}= \delta_{n,F} : & {\Z}[{\PP}^{1}(F)]\rightarrow ({\B}_{n-1}(F) 
\otimes F^{\times}),    \vphantom{\Big|} \\
 & [x]\mapsto \left\{ \begin{tabular}{ll}
0 & if $x= 0,1,\infty$, \\
$\{x\}_{n-1}\otimes x$ & otherwise,
\end{tabular}
\right.
\end{array}$$
where  $\{x\}_{n}$ denotes the class of  $x$ in $\B_{n}(F)$. \\
Although it is not used in the inductive definition,
let us define ${\RR}_{1}(F)$ to be the group generated by  $[\infty]$ and 
 $[x+y-xy]-[x]-[y]$, where $x, y\in F\backslash \{1\}$. Then  ${\B}_{1}(F)\iso F^{\times}$. \\
For $n\geqslant 2$, we define ${\A}_{n}(F)$ as the kernel of $\delta_{n}$ and ${\RR}_{n}(F)$ as the subgroup of  ${\Z}[{\PP}^{1}(F)]$ spanned by  
$[0]$, $[\infty]$ and the elements  $\sum n_{i}([f_{i}(0)]-[f_{i}(1)])$, where the $f_{i}$ are rational fractions in the indeterminate $T$, such that $\sum n_{i}[f_{i}]\in {\A}_{n}(F(T))$. 
Goncharov  proved the following basic 
\begin{lem}
For all $n \geqslant 2$, the group  ${\RR}_{n}(F)$ is contained in the kernel of  $\delta_{n}$. 
\end{lem}
\begin{proof} See \cite{g1}(Lemma 1.16, p.221) and also \cite{cathe-pol}(Proposition 1, p.1330).
\end{proof}
We then have a (cochain) complex, due to Goncharov \cite{g1,g3}, with the group $\B_n(F)$ put in degree 1,
$${\B}_{n}(F)\stackrel{\delta}{\to}{\B}_{n-1}(F) \otimes
F^{\times}\stackrel{\delta}{\to}{\B}_{n-2}(F) \otimes
\textstyle{\bigwedge}^{2}F^{\times}\stackrel{\delta}{\to}\cdots
\stackrel{\delta}{\to}{\B}_{2}(F) \otimes
\bigwedge^{n-2}F^{\times}\stackrel{\delta}{\to}
\frac{\bigwedge^{n} F^{\times}}{(2-torsion)},$$ 
with
$$\delta(\{x\}_{n-i}\otimes
y_{1}\wedge...\wedge y_{i})= \{x\}_{n-i-1}\otimes x\wedge y_{1}\wedge...\wedge
y_{i}\,,\qquad i=0,\dots,n-3,$$
and
$$\delta(\{x\}_{2}\otimes y_{1}\wedge...\wedge y_{n-2})=
(1-x)\wedge x\wedge y_{1}\wedge...\wedge y_{n-2}.$$ 
\\ \
Zagier's {\sl higher Bloch groups} \cite{z1} arise in this context as the first cohomology group of the above complex, namely
$$B_n(F)=\frac{\A_n(F)}{\CR_n(F)}\,.$$
Note the typographical difference: one has $B_n(F)\subset\B_n(F)$.
(There are in the literature several similar definitions of the ``set of relations'' $\CR_n(F)$,
denoted also $\CC_n(F)$ in \cite{z1}.) 


\begin{rem} \begin{enumerate}
\item According to Zagier's main conjecture, the groups $B_n(F)$ in the
case of a number field $F$ are presumably rationally isomorphic to $K_{2n -1}(F)_\Q$.
Using his complex, Goncharov was able to formulate a corresponding conjecture for {\sl any} field and involving the $\gamma$-filtration of the K-theory of $F$.\\
\item One of the major achievements concerning the above complexes was Goncharov's 
proof \cite{g1}  of Zagier's conjecture for $n=3$ in the course of which
he has given an explicit set of relations for (some version of) $\CR_3(F)$ which enabled him
to relate  $B_3(F)$ to (some graded piece of) the algebraic $K$-group $\,K_5(F)\,$. It is not known, however,
whether his relations generate {\sl all} functional equations for the 3-logarithm.\\
\end{enumerate}
\end{rem}

\newsubsubsection{Functions on the polylogarithmic groups}
The following proposition relates functional equations for polylogarithms and 
relations in $\B_n(F)$. (It is essentially the content of \cite{z1}, Prop.3, in the form given in \cite{g1}.)

\begin{crit} The function $\pol_n$ vanishes on $\CR_n(F)$, assuming that $F\subset \C$. \end{crit}
\
Let us end this section with a characterization of functions which actually can be defined on the corresponding $\B_n(F)$.
For $\,n\leqslant 3\,$ one knows from work of Bloch and Goncharov, respectively, a 
characterization of the measurable functions which are defined on $\,\B_n(\C)\,$:
\begin{prop} {\rm (Characterization of $\,\pol_1\,$, $\,\pol_2\,$ and $\,\pol_3\,$)}
\begin{enumerate}
\item The function $\pol_1(z)=-\log|1-z|\,$ is (up to a constant 
factor) the only measurable function  defined on $\,\B_1(\C)\,$.\\
\item  The function  $\pol_2$ is (up to a constant 
factor) the only measurable function $:\C\to\R$  which vanishes on 
$\CR_2(\C)$ and thus defines a morphism on $\,\B_2(\C)\,$.\\
\item  The space of measurable functions $:\C\to\R$  which vanish on 
$\CR_3(\C)$ and thus define a morphism on   $\,\B_3(\C)\,$, 
is two-dimensional, spanned by $\ \pol_3\ $ and $z\mapsto \log|z|\,\pol_2(z)$.
\end{enumerate}
\end{prop}
\begin{proof}
1. is classical, \ 2. has been proved by Bloch \cite{bloch1}, and 3. was given by Goncharov \cite{g1}.
\end{proof} 
%
%
\subsection{The infinitesimal polylogarithmic groups}
Cathelineau \cite{cathe-pol} has given analogues of the Goncharov complexes for infinitesimal polylogarithms whose 
cohomology is expected to be computed by some graded piece of Hochschild homology (the latter can be viewed in a sense as arising from applying 
a certain tangent functor to algebraic $K$-theory).

One defines the group $\,\beta_2(F)\,$, for $F$ any infinite field, as follows
$$\beta_2(F)=\frac{F[\Fdotdot]}{r_2(F)},$$
where $\,r_2(F)\,$ is  the kernel of the map 
$$F[\Fdotdot] \fl F^+ \otimes \U{F},\quad [a] \mapsto a\otimes a +
(1-a)\otimes (1-a).$$ If $\,\pol_2\,$ denotes the Bloch-Wigner dilogarithm
function, as defined in (\ref{classical-poly}), and if $\,F\subset \C\,$, then $\,\wt{d{\pol}_2}\,$, a somewhat modified differential
defined below, is zero on $\,r_2(F)\,$.\\
For $\,n\geqslant 3\,$, one defines inductively
$$\beta_{n}(F) =\frac{ F[\Fdotdot]}{r_{n}(F)},$$
where $\,r_{n}(F)\,$ is the kernel of the map 
$$\partial_{n} =\partial_{n,F} : F[\Fdotdot]\rightarrow (\beta_{n-1}(F)
\otimes F^{\times})\oplus ({\mathcal B}_{n-1}(F)\otimes F), $$
$$[a] \mapsto \langle a\rangle_{n-1}\otimes a + \{a\}_{n-1}\otimes (1-a),$$
and where $\,\langle a\rangle_k\,$ and $\{a\}_k$ denotes the class of $\,[a]\,$ in $\,\beta_{k}(F)\,$ and $\B_k(F)$, 
respectively.
\\

The $\,F$-vector spaces $\,\beta_{n}(F)\,$ can be viewed as  infinitesimal analogues of the 
groups  $\,{\mathcal B}_{n}(F)\,$.
The previous definition still makes sense in the case of a 
finite field $F$, but it would give $\beta_2(F)=0$. But there is also a presentation of $\beta_2(F)$ in terms of generators and relations given in \cite{cathe-sl2}(section 1, pp.52-53). As we are mainly interested by the structural properties of infinitesimal polylogarithms, we introduce the following group
\begin{df} \label{defbeta}
Let $F$ be an arbitrary field. The group $\ourbeta_2(F)$ is defined as the $F$-vector space generated by symbols $\langle a\rangle$, $a \in \Fdotdot$, subject to the relation
$$\langle a\rangle -\langle b\rangle +a\Big\langle \frac b a\Big\rangle +(1-a)\Big\langle \frac{1-b}{1-a}\Big\rangle =0,$$
for $a\neq b$.
\end{df}
We should notice that we always have a natural map $\ourbeta_2(F) \to \beta_2(F)$. 
In characteristic 0, using  \cite{cathe-pol}(section 4.2, pp. 1336-1337), we have 
\begin{prop} If $F$ is a field of characteristic 0 then the groups $\ourbeta_2(F)$ and $\beta_2(F)$ are isomorphic. 
\end{prop}


\begin{rem} \begin{enumerate} 
\item Definition \ref{defbeta} makes sense for any field.\\
\item It is 
not obvious that for a finite field $\F$ 
of characteristic $p$  we have $\ourbeta_2(\F)\neq 0$. 
It will be proven later that this is actually the case. As a counterpoint, if $F$ is a finite field or, more generally, a perfect field of characteristic $p\neq 2$, we then have $\beta_2(F)=0$ (see \cite{cathe-sl2}(Th\'eor\`eme 1, p.57)).\\
\item The infinitesimal analogue (in the sense of Cathelineau) of the above higher Bloch group $B_n(F)$ would be $\ker \partial_n/r_n(F)$
which turns out to be 0 for $n=2,3$ if $F$ is any field of characteristic 0. In fact we can show that the analogue of the  Bloch group $B_2(F)$ is given by the second Harrison homology group \cite{pev-ktheory}, proving that it is zero for any smooth $\Q$-algebra. The results and problems described  in \cite{pev-ktheory,jlc-essen}, illustrate the (presumably) close connection between infinitesimal Bloch groups and smoothness properties.
\end{enumerate}
\end{rem}

\begin{obs}\label{obs1} (Possible extension of generators in characteristic 0)
\begin{enumerate}
\item If we allow the symbols $\,\langle1\rangle_n\,$ and $\,\langle0\rangle_n\,$ in $\,\beta_n(F)\,$ then,
using the distribution relation (\ref{distri}) below, we necessarily have  $\,\langle1\rangle_n=\langle0\rangle_n=0\,$ if $n=2,3$.
\item We have $\,\langle-1\rangle_{2k+1}=0\,$ by the inversion relation.
\end{enumerate}
\end{obs}

\newsubsubsection{Functions on infinitesimal polylogarithmic groups}

The following proposition from \cite{cathe-pol} relates, for $\,F=\C\,$, functional equations for 
the infinitesimal polylogarithms and relations in the corresponding groups. 
\begin{prop}\cite{cathe-pol} For  $\,n\geqslant 2\,$, the morphism of  $\,\R$-vector spaces
$$\begin{array}{cc}
 \widehat{d{\pol }_{n}} : &\C [\C^{\scriptscriptstyle\bullet\bullet}] \longrightarrow 
\R(n-1) \\
 & b[a] \mapsto d{\pol }_{n}(a)(a(1-a)b),
\end{array}$$
is zero on $\,r_{n}(\C)\,$, hence we get a morphism
$$\wt{d{\pol }_{n}}:\beta_{n}(\C)\longrightarrow \R(n-1).$$ 
\end{prop}
\begin{rem} The definition is to be understood as follows: consider $\C$ as a 2-dimensional $\R$-vector
space with basis $(1,i)$ and with multiplication induced by the one in $\C$. Then $\pol_n$ is seen as a map from $\R^2 \to \R$, $d\pol_n(a)$ is given by the Jacobian matrix in $a$ (i.e. a row matrix of length 2). Identifying $a(1-a)b$ as a column vector relative to the basis $(1,i)$, the expression $d{\pol }_{n}(a)(a(1-a)b)$ is just the evaluation of the linear map $d\pol_n(a)$ in $a(1-a)b$ (i.e. the product of a row matrix of length 2 by a column vector of same size).

\end{rem}

\def \phil{\,G}
\begin{prop} {\rm (Characterization of $\,d\pol_2\,$)}\\
The function $d\pol_2$, restricted to $\R$, is (up to a constant factor) the only continuous function $\phil:\R^\doubledot\to \R$ 
which satisfies the equation
$$a(1-a) \phil(a) - b(1-b)\phil(b)+\frac{b(a-b)}{a} \phil\left(\frac{b}{a}\right) +\frac{(1-b)(a-b)}{1-a}\phil\left(\frac
{1-b}{1-a}\right)=0\,.$$
whenever the terms are defined.
\end{prop}
\begin{proof} Define $H(a)=a(1-a)G(a)$, $a\in\R^\doubledot$, and $H(0)=H(1)=0$, then the above functional equation is reduced to
the equation from \ref{defbeta}, for which it is well-known (cf.~e.g.~\cite{kontse}) that there is only
the differentiable function  $H(x)=-x\log|x|-(1-x)\log|1-x|$ (up to a multiplicative constant) which 
satisfies the latter equation. \end{proof}

\begin{rem} Aczel and Dhombres \cite{ad}(section 5.4, pp.66-69) have shown that if $g$ is a real  function locally integrable on $]0,1[$ and if, moreover, $g$ fulfills the  {\sl Fundamental Equation of Information Theory}, namely
$$g(x)+(1-x)g\left(\frac{y}{1-x}\right)-g(y)-(1-y)g\left(\frac{x}{1-y}\right)=0,$$
 then there exists $c\in \R$ such that $g=cH$, where $H : ]0,1[\to \R$, is the function $H(x)=-x \log(x)-(1-x) \log(1-x)$. For more detail on this topic see \cite{aivpt}.
\end{rem}

\section{Functional equations}


\begin{df}\label{def_of_fe}
A {\rm functional equation of the $n$-logarithm} resp. {\rm infinitesimal $n$-logarithm over the field $F$} is 
an element in $\CR_n(F)$ resp. in $r_n(F)$ (cf.~s.\ref{goncomplex}).

Let $F=K(t_1,\dots,t_r)$  and $K'$ be an extension of $K$. We will say that $t_1=z_1$,\dots, $t_r=z_r$, with $z_i\in K'$, is an admissible $K'$-specialisation for a functional equation $\xi(t_1,\dots,t_r)\in\CR_n(F)$ (resp. $r_n(F)$), if $\xi(z_1,\dots,z_r)$ is well defined as an element of $\ker(\delta_{n,K'})$ (resp. $\ker(\partial_{n,K'})$).
\end{df}

\begin{rem}
The restriction in the definition of a functional equation for the $n$-logarithm to {\sl rational} arguments (in the definition of $\CR_n(F)$), as opposed to algebraic arguments, is probably
not a serious one, since the corresponding polylogarithmic groups are expected to be rationally isomorphic 
(cf. e.g. \cite{g1}, pp.225, Conjecture~1.20). The above definition has the advantage of being more directly accessible to
calculations.  \\
\end{rem}

\subsection{Functional equations for classical polylogarithms}
We first list the equations which are true for general $n$: the inversion and distribution relations.
\begin{prop} {\rm (Functional equations for $\,\pol_n\,$, any $\,n\,$)}
\begin{enumerate}
\item The inversion formula:
$$\Big\{ \frac1a \Big\}_n = (-1)^{n-1}\{ a\}_n \,.$$
\item The distribution formula  
$${\left\{ a^m\right\}_n}=m^{n-1}\sum_{\zeta^m=1} 
 \,{\left\{ \zeta a\right\}}_n\, $$
holds in $\B_n(\C)$ for $\,m\in \Z\,$ and reduces to the inversion relation for $m=-1$.
\end{enumerate}
\end{prop}

\begin{rem} There is another symmetry coming from the complex conjugation:
$$\pol_n(\overline{z})=(-1)^{n-1}\pol_n(z)\,.$$
Note that this does not come from a functional equation in the above sense, since the corresponding relation  $\ \{\overline{z}\}_n
+(-1)^n\{z\}_n\ $ is not zero in $\B_n(\C)$.
\end{rem}

\newsubsubsection{The case $\,n=2$} The following functional equations are well-known for
the dilogarithm: apart from the distribution relations above it satisfies a 2-term relation relating
the arguments $x$ and $1-x$, while the most important relation (which actually characterizes $\pol_2$)
is the five term relation which allows a formulation as a 3-cocycle relation.
\begin{prop} {\rm (Functional equations for $\,\pol_2\,$)}
\begin{enumerate}
\item A two term relation. 
\begin{equation} \{ x\}_2 = -\{ 1-x\}_2\,.\end{equation}

\item The five term relation. We give two different formulations:
\begin{enumerate}
\item (as a cocycle relation in five variables): denote $\ \crr(a,b,c,d)=\frac{a-c}{ a-d}\frac{b-d}{b-c}\
$. Then 
\begin{equation}\label{five5}
 \sum_{i=1}^5 \,\,(-1)^i\,{\big\{
  \crr(x_1,\dots,\hat{x_i},\dots,x_5)\big\}_2}\,=0\,.
\end{equation}
\item in two variables (using the arguments as in Suslin's definition of the 
Bloch group; this equation is a specialization of \mbox{(a)}, putting $(x_1,\dots,x_5)=(\infty,0,1,a,b)$):
\begin{equation}\label{five2}
\left\{ a \right\}_2 \,-\, \left\{ b \right\}_2 \,+\,
\left\{ \frac{b}{a} \right\}_2 \,-\,\left\{
\frac{1-b}{1-a} \right\}_2 \,+\,\left\{
\frac{1-b^{-1}}{1-a^{-1}} \right\}_2 \,=\,0\,.
\end{equation}
\end{enumerate}
\end{enumerate}
\end{prop}

\newsubsubsection{The case $\,n=3$} For the trilogarithm one has, in addition 
to the inversion and distribution relations, an equation with 3(+1) terms
(in one variable), the well-known Kummer-Spence equation with 9(+1) terms 
(in two variables) and, most important, Goncharov's equation with 22(+1) terms
(in three variables, the ``+1'' referring to some constant term).

\begin{prop} \label{D3} {\rm (Functional equations for $\,\pol_3$)}
\begin{enumerate}
\item There is a 3-term relation 
\begin{equation}
\{ 1-x\}_3+\{ x\}_3 +\left\{
  1-\frac{1}{x}\right\}_3=\{1\}_3\,.
\end{equation}
\item The Kummer-Spence equation:  
\begin{multline}\label{kummer1}
\left\{{\displaystyle \frac {a\,(1-b)}{ b\,( 1- 
a )}} \right\}_3   +  \,\left\{{\displaystyle \frac {( 1- a )\,a}{b\,(1
 - b)}} \right\}_3 +  \,\left\{{\displaystyle \frac {a\,b}{(1 - b)\,( 1- 
a )}} \right\}_3     \\
  \mbox{}  -2 \,\left\{{\displaystyle \frac { 1- a}{1 - b}} \right\}_3  -2 \,\left\{{\displaystyle \frac {b}{b - 1}} \right\}_3 -2 \,\left\{{\displaystyle 
\frac {a}{ a- 1 }} \right\}_3  \\
  \mbox{} -2 \,\left\{{\displaystyle \frac {b}{a}} \right\}_3  -2 \,{\left\{  {\displaystyle \frac {a}{1-b}} \right\}_3 }  -2 \,{
\left\{  {\displaystyle \frac {1- a}{b}} \right\}_3 }
 +2\left\{ 1\right\}_3 =0\,.
\end{multline}
An equivalent version is given by
\begin{multline}\label{kummer2}
\left\{\frac{ x(1-y)^2}{ y(1-x)^2}\right\}_3+\left\{ x y\right\}_3+ \left\{ \frac{x}{y}\right\}_3\\- 2\left\{ \frac{y(1-x)}{y-1}\right\}_3
-2 \left\{ \frac{1-x}{1-y}\right\}_3 -2 \left\{ \frac{y(1-x)}{x(1-y)}\right\}_3 \\
-2\left\{ \frac{x-1}{x(1-y)}\right\}_3 -2\left\{ x \right\}_3 -2\left\{ y\right\}_3 
+2\left\{ 1\right\}_3 =0\,.
\end{multline}
\item Goncharov's equation: 
Set
\begin{align}\label{goncha}
f(a,b,c)=&\{a\}_3+ \left\{  \frac {b( 1- a )}{b - 1}\right\}_3
+\left\{  \frac {a( 1- b)}{ a- 1 }\right\}_3+  \left\{\frac { 1- a}{1 - abc}\right\}_3+\left\{ \frac {cb(1 - a)}{1 - abc}\right\}_3 \\
& - 
 \{ab\}_3- \left\{-  \frac {a( 1- c)( 1 - c)}{( 1- a)(1 - abc)}\right\}_3\,.
\end{align}
Then  
$$f(a,b,c)+f(b,c,a)+f(c,a,b)+\{abc\}_3 = 3 \{1\}_3\,.$$

\end{enumerate}
\end{prop}

\newsubsubsection{The case $\,n>3$} For general $\,n\,$, there are only the inversion relation
and the distribution relations known (they are the so-called trivial ones), while the existence of
non-trivial equations has only been established up to $\,n\leqslant 7\,$ (cf. \cite{gangl1}).

\subsection{Functional equations for infinitesimal polylogarithms}
Most of the functional equations for $d\pol_n$ stated in this section can be viewed as analogues of equations 
for the corresponding $\pol_n$. The main example which cannot be interpreted in this way (so far) is Cathelineau's 
equation for $d\pol_3$.\\
We first list the equations which are true for general $n$: the analogues of the  inversion and distribution 
relations.
\begin{prop} {\rm (Functional equations for $\,d\pol_n\,$, any $\,n\,$)}
\begin{enumerate}
\item The inversion formula 
\begin{equation}
\label{inver}
 a\Big\langle \frac1a \Big\rangle_n = (-1)^{n-1}\langle a\rangle_n\,.
\end{equation}
\item The distribution formula 
\begin{equation}
\label{distri}
{\left\langle a^m\right\rangle_n}=m^{n-2}\sum_{\zeta^m=1} \frac{1-a^m}{1-\zeta
  a}  \,{\left\langle \zeta a\right\rangle}_n\, 
\end{equation}
holds in $\,\beta_n(\C)\,$ for $\,m\in \Z\,$ and reduces to the inversion relation for $m=-1$. When $m=2$, we 
call this equation the {\bfseries duplication formula}.
\end{enumerate}
\end{prop}
%
%
\newsubsubsection{The case $\,n=2$} The following functional equations are true for
the infinitesimal dilogarithm: 
\begin{prop} {\rm (Functional equations for $\,d\pol_2\,$)}
\begin{enumerate}
\item The {\bfseries 2-term relation}. 
\begin{equation}  \langle x\rangle_2 = \langle 1-x\rangle_2\,.\label{twoterm}\end{equation}
\item A six term relation. Let $\,s \in F\,$. Then
\begin{equation}\label{fund1}
 (1-y)\big\langle\frac{x-s}{1-y}\big\rangle_2
+y\big\langle\frac{s}{y}\big\rangle_2 +\langle y\rangle_2
\end{equation}
is symmetric in $\,x\,$ and $\,y\,$. Specifically, we have for $s=0$ the\\
 {\bfseries fundamental equation of information theory}
\begin{equation}\label{fund3}
 (1-y)\big\langle\frac{x}{1-y}\big\rangle_2  -\langle x\rangle_2
=(1-x)\big\langle\frac{y}{1-x}\big\rangle_2  -\langle y\rangle_2
\end{equation}
which is equivalent to Cathelineau's version
\begin{equation}\label{fund2}
\langle a\rangle_2 -\langle b \rangle_2\,+\,a\Big\langle\frac b
a\Big\rangle_2
\,+\,(1-a)\Big\langle\frac{1-b}{1-a}\Big\rangle_2\,=\,0\,.
\end{equation}

\item A family of five term relations is given by taking linear combinations of the following two equations 
in five variables: 
denote $\ \crr(a,b,c,d)=\frac{a-c}{ a-d}\frac{b-d}{b-c}\ $
and $\ {\rm denom}(a,b,c,d)=(a-d)(b-c)\ $. Then  one has
\begin{equation}\label{fivex5}
 \sum_{i=1}^5 \,\,(-1)^i\,{{\rm denom}(x_1,\dots,\widehat{x_i},\dots,x_5)\,\big\langle
  \crr(x_1,\dots,\widehat{x_i},\dots,x_5)\big\rangle_2}\,=0\,,
\end{equation}

and
\begin{equation}\label{fivex5bis}
 \sum_{i=1}^5 \,\,(-1)^i\, x_i\,{{\rm denom}(x_1,\dots,\widehat{x_i},\dots,x_5)\,\big\langle
  \crr(x_1,\dots,\widehat{x_i},\dots,x_5)\big\rangle_2}\,=0\,.
\end{equation}
\item The same family of five term relations can be stated with less parameters in the arguments:
\begin{equation}
 (b+t)\langle a\rangle_2 -(a+t)\langle b\rangle_2 +(1+t)a\Big\langle\frac b
a\Big\rangle_2
\,+\,t(1-a)\Big\langle\frac{1-b}{1-a}\Big\rangle_2\,+\,b(1-a)\Big\langle\frac{a(1-b)}{b(1-a)}\Big\rangle_2\,.
\end{equation}
\end{enumerate}
\end{prop}
\begin{proof} It is a straightforward matter to check that the above elements 
lie in the kernel of $\partial_2$.  Nevertheless, we give some 
interrelationships between the various equations.\\
\begin{enumerate}
\item The symmetry of equation (\ref{fund1}) is equivalent to (\ref{fund2}):\\
We have to write the following relation
$$(1-y)\big\langle\frac{x-s}{1-y}\big\rangle_2 +y\big\langle\frac{s}{y}\big\rangle_2 +\langle y\rangle_2=(1-x)\big\langle\frac{y-s}{1-x}\big\rangle_2 +x\big\langle\frac{s}{x}\big\rangle_2 +\langle x\rangle_2\, $$
as a sum of $4$-term relations.\\
On the left hand side of the equation we add the $4$-term relation 
in the following form
\begin{align*}
-y\left\langle\frac{s}{y}\right\rangle_2-\langle y\rangle_2 &+ \langle s\rangle_2+(1-s) \left\langle\frac{1-y}{1-s}\right\rangle_2=0 \, ,
\end{align*}
and we do the same on the right hand side with $y$ replaced by $x$. 
This leaves us with another form of the 4-term relation
$$ (1-y)\left\langle\frac {x-s}{1-y}\right\rangle_2 + (1-s) \left\langle\frac{1-y}{1-s}\right\rangle_2=
(1-x)\left\langle\frac {y-s}{1-x}\right\rangle_2 + (1-s) \left\langle\frac{1-x}{1-s}\right\rangle_2$$
(to see this we should replace, in (\ref{fund2}),   $x$ by $\frac{x-s}{1-s}$ and $y$ by $\frac{y-s}{1-s}$ and use (\ref{twoterm})), thereby proving the first claim.\\
The equivalence of (\ref{fund2}) and (\ref{fund3}) is easily shown using the inversion and the
2-term relation.
\item The second family of five term relations is almost direct to deduce: the combination 
given is the sum of $t$ times the 4-term relation (\ref{fund2}) and its following equivalent 
formulation
\begin{equation} \label{fund4}
 b\langle a\rangle_2 -a\langle b\rangle_2 +a\Big\langle\frac b
a\Big\rangle_2\,+\,b(1-a)\Big\langle\frac{a(1-b)}{b(1-a)}\Big\rangle_2\,.
\end{equation}
(replace in (\ref{fund2}) $a$ and $b$ by their inverses, respectively, then multiply the result 
by $-ab$ and finally use the inversion relation on three of the ensuing terms). \\

From this, we get a very simple proof of the five term relations in cocycle form, i.e. 
(\ref{fivex5}) and (\ref{fivex5bis}): 
in each of the two versions (\ref{fund2}) and (\ref{fund4}) of the 4-term relation we put 
$a=\crr(x_1,x_2,x_3,x_4)$ and $b=\crr(x_1,x_2,x_3,x_5)$. Introducing for the moment the notation
$(ijkl):=\crr(x_i,x_j,x_k,x_l)$, we can rewrite the two equations in a concise way:
\begin{align*}
\big\langle (1234)\big\rangle &\phantom{(1234)} -\big\langle (1235)\big\rangle &+(1234)\big\langle (1245)\big\rangle  & +(1324)\big\langle (1345)\big\rangle\,, \\
(1235)\big\langle (1234)\big\rangle &-(1234)\big\langle (1235)\big\rangle &+(1234)\big\langle (1245)\big\rangle & \qquad+(1235)(1324)\big\langle (2345)\big\rangle\,.
\end{align*}
Given $\lambda\in \Z$, there is a linear combination of the two equations such that the 
coefficient of $\Big\langle\crr(x_1,x_3,x_4,x_5) \Big\rangle$ (which only occurs in the first equation)
and of $\Big\langle \crr(x_2,x_3,x_4,x_5)\Big\rangle$ (only occurring in the second equation) is
$-x_2^\lambda (x_1-x_5)(x_3-x_4)$ and $x_1^\lambda (x_2-x_5)(x_3-x_4)$, respectively.
If, for $\lambda=0$ and $\lambda=1$, we compute the coefficients of the other three arguments,
we obtain
exactly the expressions given in the claim.\\
For example, let us compute the coefficient of the first argument in the case $\lambda=1$:
the first equation is multiplied by 
$$-x_2(x_1-x_5)(x_3-x_4)\frac{(x_1-x_4)(x_3-x_2)}{(x_1-x_2)(x_3-x_4)}\,, $$
the second by 
$$x_1(x_2-x_5)(x_3-x_4)\frac{(x_1-x_5)(x_2-x_3)}{(x_1-x_3)(x_2-x_5)}\frac{(x_1-x_4)(x_3-x_2)}{(x_1-x_2)(x_3-x_4)}\,,$$
so the coefficient becomes 
$$-x_2(x_1-x_5)\frac{(x_1-x_4)}{(x_1-x_2)}(x_3-x_2)+x_1(x_2-x_5)\frac{(x_1-x_4)}{(x_1-x_2)}(x_3-x_2)\,,$$ 
which is equal to $\ x_5(x_1-x_4)(x_2-x_3)\,.$
\end{enumerate}
\end{proof}

\begin{rem}
The generalized version  of the fundamental equation of information theory, namely (\ref{fund1}), is
equivalent to the one given by both Kontsevich and Cathelineau (referring to
Acz\'el-Dhombres), as was shown in the proof (part 1.) above. At first glance, it is somewhat surprising 
that we do not get anything new although we can achieve to insert a third (non-homogenizing) parameter---but
there are related phenomena known for the five term relation. 
In particular, we do not gain new information for information theory.
\end{rem}

\newsubsubsection{The case $\,n=3$} For the infinitesimal trilogarithm one has an equation with three terms
(in one variable), a ``derived version'' of the Kummer-Spence equation with 
eight terms 
(in two variables) and, most important, Cathelineau's equation with 22 terms
(in three variables).\\
 
The proposition below gives complementary information on $\beta_3(F)$.
\begin{prop} \label{dD3} {\rm (Functional equations for $\,d\pol_3$)}
\begin{enumerate}
\item There is a 3-term relation 
\begin{equation}\label{three}
\langle 1-x\rangle_3-\langle x\rangle_3 +x\left\langle
  1-\frac{1}{x}\right\rangle_3=0\,.
\end{equation}
\item The Kummer-Spence analogue: the $\,F$-linear combination 
\begin{multline}\label{kummerana1}
  \,{\displaystyle \frac {(1 - b)\,b}{1- b - a}
\,\left\langle{\displaystyle \frac {( 1- a )\,a}{b\,(1
 - b)}} \right\rangle_3 } +  \,{\displaystyle \frac {(1 - b)\,(1 - a
 )}{1- b - a}\,\left\langle{\displaystyle \frac {a\,b}{(1 - b)\,( 1- 
a )}} \right\rangle_3 }  \\
  \mbox{}  + \,(1-b)\,\left\langle{\displaystyle \frac { 1- a}{1 - b}} \right\rangle_3 
 - \,(1 - b)\,\left\langle{\displaystyle \frac {b}{b - 1}} \right\rangle_3 - \,(1 - a)\,\left\langle{\displaystyle 
\frac {a}{ a- 1 }} \right\rangle_3   \\
  \mbox{} 
  - \,a\,\left\langle{\displaystyle \frac {b}{a}} \right\rangle_3 + \,{\displaystyle \frac {(a
 - b - 1)\,(1 - b)}{1-b-a}\,\left\langle  {\displaystyle \frac {a}{1-b}} \right\rangle_3 } 
- \,{\displaystyle \frac {(a  - b+ 1)\,b}{ 1- b - a}\,
\left\langle  {\displaystyle \frac {1- a}{b}} \right\rangle_3 }
\end{multline}
vanishes in $\,\beta_3(F)\,$. An equivalent version, denoted $KS(x,y)$, is given by
\begin{multline}\label{kummerana2}
\left\langle x y\right\rangle_3+y \left\langle \frac{x}{y}\right\rangle_3-(1-y) \left\langle \frac{y(1-x)}{y-1}\right\rangle_3\\
+(1-y) \left\langle \frac{1-x}{1-y}\right\rangle_3 - x(1-y) \left\langle \frac{y(1-x)}{x(1-y)}\right\rangle_3 \\
+x(1-y)\left\langle \frac{x-1}{x(1-y)}\right\rangle_3 -(1+y)\left\langle x \right\rangle_3 -(1+x)\left\langle y\right\rangle_3 
\end{multline}
\end{enumerate}
\end{prop}
\begin{proof} The 3-term equation and (\ref{kummer2}) will follow directly
from the equation in the next  proposition. The equivalence of the two Kummer-Spence analogues becomes
  evident after applying the change of variables $\,x=\frac{a}{1-b}\,$, 
$\,y=\frac{1-a}b\,$, and multiplying the result by $\,\frac{b(1-b)}{1-a-b}\,$.
\end{proof}

We can also notice the following formal property, that we will give as
\begin{lem}\label{3-term_implies_inversion} In $\beta_3(F)$, the inversion formula is a consequence of the 3-term equation.
\end{lem}
\begin{proof} Add the 3-term equation to its variant where $x$ is replaced by $1-x$. Four of
the terms cancel and the remaining two give the inversion relation.
\end{proof}
Cathelineau has given a 22-term equation which completely describes the
set of relations for the infinitesimal polylogarithmic group $\,\beta_3(F)\,$: 
In order to state it conveniently, we use his notation for a distinguished linear combination of seven terms
\begin{equation}\label{doublebra}
[[a,b]]=(b-a)\tau(a,b) + \frac{1-b}{1-a}\sigma(a)+\frac{1-a}{1-b}\sigma(b)\,,
\end{equation}
where we have set
$$\tau(a,b)=\frac{[a]}{1-a}-\frac{[b]}{1-b}+ \frac{a}{a-b}\Big[\frac{ b}{ a}\Big]  -\frac{1-a}{b-a}\Big[\frac{1-b}{1-a}
\Big] + \frac{b(1-a)}{b-a}\Big[\frac{a(1-b)}{b(1-a)}\Big]\,,$$
($\tau$ arises by taking the five term relation (\ref{five2}) and multiplying each $[z_i]$ with the coefficient $\frac{1}{1-z_i}$) and
$$\sigma(a)=a[a]+(1-a)[1-a]\,.$$
Then we can state the 22-term relation as follows

\begin{df} We define the formal expression $J(a,b,c)$ in the  indeterminates $a,b,c$ as
\[ J(a,b,c)=[[a,c]]-[[b,c]]+a\Big[\Big[\frac{b}{a},c\Big]\Big]+(1-a)\Big[\Big[\frac{1-b}{1-a},c\Big]\Big]\, .\]
\end{df}
\begin{rem}
\begin{enumerate}
\item Writing out all the terms, we obtain 22 different arguments: 
\begin{align*}
J(a,b,c)
=\ &\ c\left[ a\right] - 
c\left[ b\right] +(a-b+1)\left[ c\right]\\[0.1cm]
&+(1-c)\left[ 1-a\right] -(1-c)\left[ 1-b\right] 
+ (b-a)\left[ 1-c\right]\\[0.1cm]
&-a\left[ \frac{c}{a}\right]  
+b\left[ \frac{c}{b}\right] 
+ca\left[ \frac{b}{a}\right]\\[0.1cm]
&-(1-a)\left[ \frac{1-c}{1-a}\right]
+(1-b)\left[ \frac{1-c}{1-b}\right] 
+c(1-a)\left[ \frac{1-b}{1-a}\right]\\[0.1cm]
&+c(1-a)\left[ \frac{a(1-c)}{c(1-a)}\right] 
-c(1-b)\left[ \frac{b(1-c)}{c(1-b)}\right]\\[0.1cm]
&-b\left[ \frac{ca}{b}\right] 
-(1-b)\left[ \frac{c(1-a)}{1-b}\right] \\[0.1cm]
&+(1-c)a\left[ \frac{a-b}{a}\right] 
+(1-c)(1-a)\left[ \frac{b-a}{1-a}\right] \\[0.1cm]
&-(a-b)\left[ \frac{(1-c)a}{a-b}\right] 
-(b-a)\left[ \frac{(1-c)(1-a)}{b-a}\right] \\[0.1cm]
&+c(a-b)\left[ \frac{(1-c)b}{c(a-b)}\right] 
+c(b-a)\left[ \frac{(1-c)(1-b)}{c(b-a)}\right].
\end{align*}
\item When $a,b,c$ are elements of an arbitrary  field $F$, we will still use the notation $J(a,b,c)$ for the evaluation of $J$ in the specified values.
\end{enumerate}
\end{rem}
\ 
\begin{thm} \label{main_jlc_result_for_beta3}{\rm (Cathelineau, \cite{cathe-pol}, Corollaire 1, p.1345)}\ Let $F$ be a field of characteristic zero.
\begin{enumerate}
\item The image of $\,J(a,b,c)\,$ under the projection $\ F[\Fdotdot]\to
\beta_3(F)\ $
is zero.\\
\item Furthermore, $\,J(a,b,c)\,$, together with its specializations to $c=a,$ $b,$ $\frac a b$ or 
$\frac{1-a}{1-b}$, respectively, and the inversion relation generate the set of relations which define 
$\,\beta_3(F)\,$. Here we understand $\langle 1\rangle_3=0$.
\end{enumerate}
\end{thm}
\begin{rem}\label{specialization}\begin{enumerate}
\item In the presentation of \cite{cathe-pol}, Corollaire 1, one can replace his equation 1) coming from
$[[a,b]]-[[b,a]]$ by the shorter inversion relation (\ref{inver}). (Proof: add his equation 1) to the same relation 
where $a$ and $b$ are replaced by $\frac{1}{a}$ and $\frac{1}{b}$ and where the result is multiplied by $ab$.)\\

 \item The combinations $[[a,c]]+a[[1/a,c]]$ and $[[a,c]]-[[1-a,c]]$ give versions of the Kummer-Spence 
analogue. Since, e.g., $\langle a\rangle_2-\langle 1-a\rangle_2=0$ results formally from the four term relation
(at least up to 2-torsion), we get the Kummer-Spence analogue directly from $J(a,b,c)$.\\

\item By Observation \ref{obs1}, one can introduce
elements $[a]$ for $a=0,1$ and set their image in $\beta_3(F)$ equal to zero. What is more, one can add a 
formal generator $[\infty]$ as well, for which we only require $0[\infty]=0$.
One can then formally deduce the 3-term equation (\ref{three}) by specializing $a=1$ in 
$\frac{1}{b-1}J(a,b,c)$ and one obtains the
Kummer-Spence analogue (\ref{kummer2}) by specializing $a=0$ in $(1-x)(1-y)J(a,(1-x)^{-1},(1-y)^{-1})$,
(these specializations are not allowed in Cathelineau's context, but will make sense in the ``finite polylog'' case below).\\

\item A different way to obtain the Kummer-Spence analogue is to symmetrize, i.e.~to form
$$J(a,b,c)+J(b,a,c)+c\big(J(a,b,\frac 1 c)+J(b,a,\frac 1 c)\big)\,,$$
and then to check that one obtains the difference of two Kummer-Spence
analogues $KS(c,\frac b a)-KS(c,\frac {1-b}{1-a})$.\\
\item Alternatively, one can deduce the Kummer-Spence analogue or the 3-term relation (non-explicitly) from $J(a,b,c)$
by simply checking that the corresponding linear combinations
lie in the kernel of $\,\partial_3\,$, and then use Cathelineau's theorem to
deduce that each such combination must be a consequence of $\,J(a,b,c)\,$.\\

\item In the case of the classical trilogarithm, Goncharov has given
a new functional equation in 22(+1) terms which {\sl presumably} generates all 
functional equations for $\pol_3$, i.e.~the kernel of $\delta_3$, but there are (infinitely many) functional
equations (cf. \cite{wojt3}, \cite{gangl3}) which are not known to be formal consequences of it. 
Cathelineau's result in the infinitesimal setting is stronger in the sense that it actually generates the kernel of $\partial_3$.
\end{enumerate}
\end{rem}
One of the major consequences of  Theorem \ref{main_jlc_result_for_beta3} is that it allows us to give a general definition for $\ourbeta_3$.
\begin{df}\label{defbeta3}  
Let $F$ be an arbitrary field. The group $\ourbeta_3(F)$ is defined as the $F$-vector space generated by symbols $[a]$, $a \in F$,
subject to the relations $\,J(a,b,c)\,$, together with its specializations to $c=a,$ $b,$ $\frac a b$ or 
$\frac{1-a}{1-b}$, respectively, the inversion relation and $[1]=[0]=0$.
\end{df}
If in $\beta_3(F)$ we introduce elements $[a]$ for $a=0,1$, we then have, in virtue of 
(\ref{obs1}), a surjective map of $F$-vector spaces $\ourbeta_3(F) \to \beta_3(F)$, which is an isomorphism in characteristic 0. As in the case $n=2$, if $F$ is a finite field of characteristic $p$, $\beta_3(F)=0$ but it will be shown in part III that $\ourbeta_3(F)\neq 0$. The groups $\ourbeta_n(F)$, for $n=2,3$, measure how much the group $\beta_n(F)$ deviates from being generated by the main functional equations of infinitesimal polylogarithms.

\newsubsubsection{The case $\,n>3$} For general $\,n\,$, there are only the inversion relation
and the distribution relations, as seen in   (\ref{distri}), known. For each functional equation of the
corresponding classical polylog, using the ``derivation map'' described in  the section \ref{deriving}, there is associated a functional 
equation (actually many) for the infinitesimal polylogarithm.
From what has been stated above for the classical case, this 
means that at least up to $\,n=7\,$ there are non-trivial ones.



\part{The Results}
\section{Finite versions of polylogarithms and their functional equations}\label{finite}


In this section we will study what we can call {\it finite analogs} of the polylogarithms and also the groups $\ourbeta_n(F)$ for $n=2,3$ in the case where $F$ is a field of characteristic $p\neq 2$ (eventually finite). We will show that for $n=2,3$ the finite analogs of the polylogarithms define functions on $\ourbeta_n(F)$, showing that surprisingly they behave like the infinitesimal polylogarithms. As for the previous cases, we will show, at least in low dimension, that these {\it finite polylogs} are uniquely characterized by their functional equations.

For the remainder of the paper, let us fix an odd prime $\,p\,$.
We shall work over an arbitrary field $\F$ of characteristic $p$.
\subsection{Definition and first properties of finite polylogarithms}



\begin{df} For any field $F$ of characteristic $p$, the {\rm $\,n$th finite polylogarithm} or {\rm
finite $n$-logarithm} is given by the following polynomial in $\,\F[T]\,$:
$$£_{n}(T)=\sum_{k=1}^{p-1} \frac{T^k}{k^n}\,.$$
\end{df}

\begin{nota}For the remainder of this paper, we will denote $\,\wt{ P}\,$ the function associated to the polynomial $\,P\,$.
\end{nota}

\begin{rem}\begin{enumerate}
\item   ``Extension by periodicity''\\
If $\F$ is of characteristic $p$, it has $\F_p$ as prime subfield, which is fixed by the Frobenius morphism $x \mapsto x^p$. As a result we have the $(p-1)-$periodicity  $£_{n+p-1}=£_n$, and we need only consider $n<p$. 
\item It is important to notice that the functions $\widetilde{£_n}$ are not identically zero on $\F$.
\end{enumerate}
\end{rem}
The following differential equation relates the finite polylogarithms of different
orders (just like in the classical case)
\begin{align*}
d£_{n}(U) 
&=£_{n-1}(U)\,d\log(U) ,
\end{align*}
where we denoted $\frac{dU}{U}$ by $d\log(U)$. 
Extending this formally, it is convenient to introduce the following notation:
\def \hatl {\widehat{£}}
\begin{df}\label{dlog} Let $F$ be a field of characteristic $p$.
Define the following {\rm ``Frobeniizing'' map}
\begin{align*}
\hatl_m \, : \, F[\Fdotdot] &\to F\,,\\
c[f]&\mapsto c^p\pounds_m(f)\,.
\end{align*}
\end{df}
One observes that, for any  $c$ and $f$ in $F$, the differential operator $\frac{\partial}{\partial x}$ 
acts linearly on the coefficient $c$ of $\hatl_m\big(c[f]\big)$ and, as above, like $d\log$ on the generator $[f]$:
$$\frac{\partial}{\partial x}\hatl_m\big(c[f]\big)=\hatl_{m-1}\big(c[f]\big) \frac{\partial}{\partial x}\log(f)\,.$$

\begin{obs} \begin{itemize}
\item[0.] For $n=0$ we have 
$$£_{0}(T)=\frac{T-T^p}{1-T},$$
and therefore 
\begin{equation}\label{L0}T£_{0}(1-T)=-(1-T)£_{0}(T).
\end{equation}
\item[1.]  For $n=1$, by expanding $(1-T)^p$ and noticing that  $\frac1 p
{p \choose k}=\frac 1 k{p-1\choose k-1} =\frac{(-1)^{k-1}}k$, we get  a simple (and well-known) formula 
$$ \sterling_1(T)\equiv \frac{1-T^p-(1-T)^p}{p}\pmod p\,.$$
Note that the term on the right hand side occurs in the polynomials which define the sum of 
the two Witt vectors $(1,0,\dots,0,\dots)$ and $(-T,0,\dots,0,\dots)$.
\end{itemize}
\end{obs}

\subsection{Functional equations for finite polylogarithms} {\it A priori}, there seems to be at least 
two natural candidates for functional equations for the finite $n$-logarithm: we could ask for 
linear combinations 
$\sum_i c_i [x_i]$ such that $\widetilde{£_n}$ vanishes for all specializations of the parameters 
which ``make sense'' (i.e. no term ``$\frac00$'' occurs); we will call those combinations weak
functional equations. But this definition has the disadvantage that there are too many ambiguities involved (just think of a coefficient that is
divisible by $x^p-x$). Instead, we will impose the stronger property that $\sum_i c_i{£_n}(x_i)$ vanishes
as a rational expression, and by multiplying with the common denominator, we can even assume it to
vanish as a polynomial. 

\begin{df}
A {\rm functional equation in the strong sense} for the finite $n$-logarithm over a field $F$ of characteristic $p$ is a finite
linear combination 
$\sum_i c_i \sterling_n(x_i)\in F(t)[F(t)]$ which vanishes identically as a polynomial.\\
A {\rm functional equation in the weak sense} is a finite
linear combination 
$\sum_i c_i \sterling_n(x_i)\in F(t)[F(t)]$ which vanishes for each specialization of parameters
which makes sense.
\end{df}
In the following we list a number of equations which are identical to the ones for the infinitesimal
polylogarithms, apart from ``Frobeniizing'' the coefficients (i.e.~raising them to the $p$th power). 
The proofs will be postponed to \S\ref{theproofs}.

\newsubsubsection{General functional equations for $\pounds_n$} 
\begin{prop} \label{generic}
Let $\,n\in\Z\,$ be arbitrary. We have the following identities 
\begin{enumerate}
\item Inversion formula: $\,£_{n}(T)=(-1)^{n} T^p \, £_{n}\left(\dfrac{1}{T}\right)$.\\
It can be viewed as a special case $\,(m=-1)\,$ of the following
\item Distribution formulae: assume $F$ contains a primitive $m$th root of unity. Then 
$${£}_{n}(T^m)=m^{n-1} \sum_{\zeta^m=1} \frac{1-T^{pm}}{1- \zeta^p T^p}\,{£}_{n}(\zeta T)\, .$$
\item Special values: $\wt{{£}_{n}}(1)=0\,$ if $\,(p-1)\hskip -4pt\not|n\,$and $=-1$ else, while  $\wt{{£}_{2n}}(-1)=0$ for any $n$. Let $B_j=$ be the $jth$ Bernoulli number and set $G_j=2(1-2^j)B_j$. Then  for $0<m<(p-1)$ we have  that $\,\wt{{£}_{p-m}}(-1)=\dfrac{G_{m}}{m}$.
\end{enumerate}
\end{prop}
\begin{rem} Notice that the numbers $G_m$ are integers by virtue of classical results (for instance it is a consequence of the Theorem of von Staudt-Clausen\cite{washington}(Theorem 5.10, p.56)). These numbers are called the Genocchi numbers and we have $mG_{p-1+m}=(m-1)G_m$ mod $p$ which is nothing else than the
famous Kummer congruence for Bernoulli numbers. 
\end{rem}

Still mirroring the set-up in the infinitesimal case, we now state several functional equations 
specific to $\,n=1,2\,$.


\newsubsubsection{Equations for $\,£_{1}$}
\begin{prop}\label{L1.1}
\begin{enumerate}
\item The 2-term relation: $£_{1}(T)=£_{1}(1-T)$.
\item The generalized fundamental equation of information theory: let $\,s\,$,
$\,x\,$ and $\,y\,$ be indeterminates. The expression
$$ H(x,y,s)=(1-y)^p£_{1}\left(\frac{x-s}{1-y}\right) +y^p£_{1}\left(\frac{s}{y}\right) +£_{1}(y)$$
 in $\ \F[x,y,s]\ $ is symmetric in $\,x\,$ and $\,y\,$.
Specifically, we have 
\begin{equation}\label{fundfinite}
£_{1}(a)- £_{1}(b) \,+\,a^p\,£_{1}\left(\frac b
a\right)
\,+\,(1-a)^p\,£_{1}\left(\frac{1-b}{1-a}\right)\,=\,0\,.
\end{equation}
\item The five term relations. \\
Denote $\ \crr(a,b,c,d)=\frac{a-c}{ a-d}\frac{b-d}{b-c}\ $
and $\ {\rm denom}(a,b,c,d)=(a-d)(b-c)\ $. Then  we have the polynomial identities in
$\F[x_1,\dots,x_5]$
$$ \sum_{i=1}^5 (-1)^i\,{\big({\rm denom}(x_1,\dots,\widehat{x_i},\dots,x_5)\big)^p\,£_{1}\big(\crr(x_1,\dots,\widehat{x_i},\dots,x_5)\big)}\,=0\,,$$
and
$$ \sum_{i=1}^5 (-1)^i\,x_i^p\,{\big({\rm denom}(x_1,\dots,\widehat{x_i},\dots,x_5)\big)^p\,£_{1}\big(\crr(x_1,\dots,\widehat{x_i},\dots,x_5)\big)}\,=0\,.$$
\end{enumerate}
\end{prop}
\begin{cor} The $\F$-vector space $\ourbeta_2(\F)$, as defined in (\ref{defbeta}),  is of dimension at least 1. If, moreover, $\F$ is a perfect field, then $\ourbeta_2(\F)=\F$.
\end{cor}
\begin{proof} According to Proposition \ref{L1.1}, the function $\widetilde{£_1}$ is a well-defined function on $\ourbeta_2(\F)$, and as it is not identically zero on $\F$, the dimension of $\ourbeta_2(\F)$ is non-zero. By \cite{cathe-sl2}(Th\'eor\`eme 1, p.57), we know that $\beta_2(\F)=0$. But as the relations in $\beta_2(\F)$ are given by the 4-term equation (i.e. the Fundamental Equation of Information Theory) and the relation $\sum_{k=2}^{p-1} [k 1_\F]$,
(see section 1.1 and also Sah's Lemma in \cite{cathe-sl2}(pp.52-53)), and as we further know, again by Sah (see the remark on p.53 in op. cit.), that these two relations are independent, we can conclude that the kernel of the map $\ourbeta_2(\F)\to \beta_2(\F)$ is generated by the element $\sum_{k=2}^{p-1} [k 1_\F]$. Evaluating $\widetilde{£_1}$ on this element shows that it is non-zero, which ends the proof.
\end{proof}
\newsubsubsection{Equations for $\,£_{2}$}
In this subsection we will give answers to the question raised by Kontsevich in \cite{kontse}. Notice that we need to assume $p>3$ throughout.
\begin{prop}\label{L2.0} The 3-term relation and the Kummer-Spence analogue are functional equations
for $£_2$.
\end{prop}
\begin{proof} This is a consequence of the following theorem, together with remark
(\ref{specialization}). \end{proof}
\begin{thm} \label{L2.1}
 The image of $\,J(a,b,c)\,$ under the map $\widehat{\pounds}_2$
is a polynomial which is identically zero in $\F[a,b,c]$.
\end{thm}
\begin{rem}
\begin{enumerate} \item
By \S\ref{red_mod_p}, there is a better answer to Kontsevich's question, at least ``quantitatively'': each functional equation for $d \pol_3$ induces a functional equation (in the weak sense) for $£_2$. This is true in 
particular for the 3-term equation and the Kummer-Spence analog. \\
\item One can find further equations (in the strong sense) for $£_2$ and in 
general for $£_n$ with $n\geqslant 3$ in   \cite{aivpt}.
\end{enumerate}
\end{rem}

\noindent By similar arguments as in the proof for $£_1$, we get
\begin{cor} The $\F$-vector space $\ourbeta_3(\F)$ is of dimension at least 1.
\end{cor}

\subsection{Characterization of finite polylogarithms} 
We can characterize $\pounds_1$ and $\pounds_2$ by the functional equations they satisfy.
\begin{prop} \label{unicity_sterling1} The space (over $\F$) of solutions of the ``fundamental equation of information theory'' is of dimension $1$ generated by $£_1$.
\end{prop}
\begin{proof} 
Set
  $\,f(T)=\sum_{i=0}^{p-1}a_i T^i \in \F_q[T]\,$, and suppose that $\,f\,$ verifies
  $\,f(0)=0\,$ and the following identity in $\,\F_q[x,y]\,$
$$f(x)+(1-x)^pf\left( \frac{y}{1-x}\right) -f(y)-(1-y)^pf\left( \frac{x}{1-y}\right)=0 \,.$$
Differentiating the previous equation with respect to $\,x\,$ gives,
$$df(x)+\frac{y (1-x)^p}{(1-x)^2}\, df\left( \frac{y}{1-x}\right)-\frac{(1-y)^p}{1-y}
\,\,df\left( \frac{x}{1-y}\right)=0,$$
with $\,df(T)=a_1+\sum_{i=2}^{p-1}i a_iT^{i-1}\,$ and thus $\,df(0)=a_1\,$. Setting $\,x=0\,$ in the previous identity gives
$$a_1+y\, df(y)-\frac{1-y^p}{1-y}\,a_1=0.$$
But as $\,\frac{1-y^p}{1-y}=\sum_{i=0}^{p-1}y^i\,$, the previous equality implies
$\,a_i=\frac{a_1}{i}\,.$
In other words, since $\,f(0)=0\,$ we have $\,f=a_1 £_{1}\,$, which proves the claim.
\end{proof}

In fact we have  a stronger statement
\begin{prop} \label{strong_unicity_of_sterling1} The 2-term equation, the inversion  and the duplication formulae characterize altogether $£_1$.
\end{prop}
\begin{proof} It is a consequence of the following lemma
\begin{lem} \label{sequence_for_sterling1} Suppose that $a_k$ is a sequence of integers with $k=1,\dots , p-1$ ($p$ an odd prime fixed), which fulfills the following rules
\[ a_k = \begin{cases} -\frac{1}{2} \sum_{i=k+1}^{p-1} a_i \binom{i}{k}, \quad &\textrm{if} \; k \; \text{is odd},\\[3pt]
\frac{1}{2} a_\frac{k}{2} \quad &\text{otherwise},
\end{cases}
\]
and $a_{p-k}=-a_k$ for all $k=1,\dots,p-1$. 
Then $a_k=\frac{a_1}{k}\in \F_p$ for all $k=1,\dots, p-1$.
\end{lem}
\begin{proof}[Proof of the lemma]
The proof goes by descending induction starting from $p-1$. First we notice that by the third rule, we have $a_{p-1}=-a_1=\frac{a_1}{p-1}$ modulo $p$. Suppose that $a_i=\frac{a_1}{i}$ modulo $p$ for all $i>k$. Now compute $a_k$ modulo $p$. Observe that we can assume $k\leqslant p-3$, since we can compute from the rules $a_{p-1}$ and $a_{p-2}$. If $k$ is odd then by the first rule we deduce directly $a_k$, but we still have to show that  $a_k=\frac{a_1}{k}$ modulo $p$. This is done via the
\begin{sublem}\label{thesublemma1} If $k$ is odd and $a_i=\frac{a_1}{i}$ modulo $p$ for all $i>k$. Then $a_k=\frac{a_1}{k}$ modulo $p$.
\end{sublem}
\begin{proof}[Proof of the sub-lemma]
We have to show that, modulo $p$, 
\[ \frac{a_1}{k} = -\frac{1}{2} \sum_{i=k+1}^{p-1} \frac{a_1}{i} \binom{i}{k},\]
or equivalently, assuming $a_1\neq 0$, that
\[ -2 =\sum_{i=k+1}^{p-1} \frac{k}{i} \binom{i}{k} .\]
But 
\[ \frac{k}{i} \binom{i}{k}=\binom{i-1}{k-1} .\]
Using the usual rule  $\binom{m}{n}=0$ if $n>m$, we have
\[\sum_{i=k+1}^{p-1} \frac{k}{i} \binom{i}{k}=\sum_{i=0}^{p-2}\binom{i}{k-1}-1 .\]
But $\sum_{i=0}^{p-2}\binom{i}{k-1}=\binom{p-1}{k}$, and as, modulo $p$, we have $\binom{p-1}{k}=(-1)^k$, we finally get, using the fact that $k$ is odd, the desired identity.
\end{proof}
Now return to the proof of the lemma and suppose that $k$ is even. If $k=2$ then the process ends, so we can suppose that $k>3$. The idea is to show that we can compute directly $a_{k-1}$ and to deduce $a_k$ from the first rule (we will still need to show the desired property). As $k$ is even, $k-1$ is odd and thus $p-k+1$ is even. Thus by the third rule we have $a_{k-1}=-a_{p-k+1}$ and by the second rule we have 
\[a_{p-k+1}=\frac{1}{2} a_{\frac{p-k+1}{2}} .\]
But there exists $j \in \N$ such that $p=k+j$ with $3 \leqslant j <p$ (because $k\leqslant p-3$). Hence, applying once again the third rule gives
\[ a_{\frac{p-k+1}{2}} =-a_{p-\frac{p-k+1}{2}} .\]
But 
\[p-\frac{p-k+1}{2}=\frac{p+k-1}{2}=k+\frac{j-1}{2} .\]
And as $j\geqslant 3$, we have $\frac{j-1}{2} \geqslant 1$, which means, applying the induction, that $a_ {p-\frac{p-k+1}{2}}$ is already known. We then get the value of $a_{k-1}$ and by applying the first rule to it we deduce the value of $a_k$. Now to finish the proof we need to show that, in this case $a_k=\frac{a_1}{k}$ modulo $p$. Notice  that we can also assume by the induction that $a_i=\frac{a_1}{i}$ modulo $p$ for all $i>k$. First we show that in the previous process, we get $a_{k-1}=\frac{a_1}{k-1}$ modulo $p$. Indeed, by the induction we have
\[ a_{p-\frac{p-k+1}{2}} = \frac{a_1}{p-\frac{p-k+1}{2}} .\]
Thus
\[ a_{\frac{p-k+1}{2}} = - \frac{a_1}{p-\frac{p-k+1}{2}} .\]
And finally
\begin{align*}
a_{k-1} &= -a_{p-k+1},\\
 &= \frac{1}{2} \left(\frac{a_1}{p-\frac{p-k+1}{2}}\right),\\
&= \frac{a_1}{k-1} .
\end{align*}
To conclude we need to prove a variant of Sub-Lemma \ref{thesublemma1}.
\begin{sublem}
\label{thesublemma2} Suppose that $k$ is even, $a_i=\frac{a_1}{i}$ modulo $p$ for all $i>k$ and $a_{k-1}=\frac{a_1}{k-1}$ modulo $p$. Then $a_k=\frac{a_1}{k}$ modulo $p$.
\end{sublem}
\begin{proof}[Proof of Sub-Lemma \ref{thesublemma2}]
We have the equality
\[ \frac{a_1}{k-1}= -\frac{1}{2} a_k k -\frac{1}{2} \sum_{i=k+1}^{p-1} \frac{a_1}{i} \binom{i}{k-1} .\]
Using the same arguments than in the proof of Sub-Lemma \ref{thesublemma1}, we get the following identities,
\begin{align*}  \sum_{i=k+1}^{p-1} \frac{k-1}{i} \binom{i}{k-1} &=\sum_{i=0}^{p-2} \binom{i}{k-2} -\binom{k-1}{k-2}-\binom{k-2}{k-2},\\
&=(-1)^{k-1}-(k-1)-1,\\
&=-1-k, \text{ as } k \text{ is even}.
\end{align*}
And we finally have
\[ \frac{a_1}{k-1}= -\frac{1}{2} a_k k+\frac{(1+k) a_1}{2(k-1)},\]
from which we deduce $a_k = \frac{a_1}{k}$.
\end{proof}
Hence the proof of Lemma \ref{sequence_for_sterling1} is complete.
\end{proof}
Back to the proof of Proposition \ref{strong_unicity_of_sterling1}. Suppose that $P(T)=\sum_{i=0}^{p-1}a_i T^i\in \F_p[T]$ verifies the conditions of the proposition. Then applying the three equations to $P$ gives $a_0=0$, and  the other coefficients $a_i$ fulfill  the  rules described in the Lemma \ref{sequence_for_sterling1}.
\end{proof}

\begin{rem} ``Cohomological characterization of $£_1$'' \\ 
Kontsevich showed that $\sterling_1$ gives a non-zero 2-cocycle in $H^2(\Z/p,\Z/p)$. Since the latter group
is isomorphic to $\Z/p$, this characterizes $\sterling_1$ up to a scalar.
\end{rem}

\subsection{Space of solutions for equations associated to $£_2$}
As $\,J(a,b,c)\,$ is the main relation for $\,\ourbeta_3\,$, we can expect that it characterizes $\,£_{2}\,$. In fact, 
we can first give a family of polynomials (which form a space of dimension growing linearly with $p$) 
and then characterize  $\,£_{2}\,$ by imposing also the duplication relation (i.e. the distribution relation
for $£_{2}$  with $m=2$). Since these two equations are consequences of the Kummer-Spence analogue, and the
latter in turn is a consequence of $J(a,b,c)$, we are done.
\begin{prop} The dimension of the $\F_p$-space of solutions associated to the equation
\begin{equation}T^p P\left(1-\frac{1}{T}\right)-P(T)+P(1-T)=0\label{3-term}\end{equation}
grows with $p$ and is at least of dimension $\frac{p-1}{3}+1$. The family of  polynomials
$$\tau_{i,p}(T)=T^{i}(1-T)^{i}(T^{p-3i}+(-1)^i),$$
with $i\in \N$ such that the valuation of $\tau_{i,p}$ is $\geqslant 0$ (for instance if $i\leqslant \floor{\frac{p}{2}}$), is a solution of (\ref{3-term}). Moreover,
for $i=0,\dots,\frac{p-1}{3}$, this family is  free.
\end{prop}
\begin{proof} The fact that $\tau_{i,p}$ fulfills  (\ref{3-term}) is a
direct computation.
For $i=1,\dots,\frac{p-1}{3}$, the family is free for degree reasons, since
$deg(\tau_{i,p})=p-i$. Furthermore $\tau_{0,p}$ does not belong to this family for valuation reasons.
\end{proof}
\begin{rem} 
\begin{enumerate}
\item We already know, by Lemma \ref{3-term_implies_inversion}, that the inversion formula is a consequence of the 3-term equation. But a straightforward computation shows that the polynomials $\tau_{i,p}$ fulfill the inversion formula for $£_2$.
\item In fact the rank of the family $\tau_{i,p}$ is greater than $\frac{p-1}{3}$, but the proof is a little bit more involved. We can also notice that $£_2$ is never expressible in terms of  $\tau_{i,p}$ if $i$ runs only through $0,\dots,\frac{p-1}{3}-1$.
\end{enumerate}
\end{rem}

\noindent Thus the 3-term equation is insufficient for the characterization of $£_2$. Nevertheless, we have the following main result

\begin{thm}\label{unicity_of_sterling2} Let $P$ be a polynomial of $\F[T]$  of degree less than or equal to $p-1$. Set $h=TP'$. Then if $P$ fulfills the  duplication relation and the 3-term equation, and if moreover $h$ fulfills the 2-term equation then $P$ is equal, up to a multiplicative constant, to $£_2$.
\end{thm}
\begin{proof} Let $P$ be a polynomial of degree $\leqslant p-1$, and suppose that $P$ fulfills the following two equations
\begin{align}
\label{EQ1} T^p P\left(1-\frac{1}{T}\right)-P(T)+P(1-T)&=0,\\
\label{EQ2} 2 (1+T^p) P(T)+2 (1-T^p) P(-T) - P(T^2) &=0 .
\end{align}
Then observing that we can deduce the inversion formula as a consequence of the 3-term, and taking the derivative with respect to these equations shows that $h$ fulfills the inversion formula and the duplication formula. As, by hypothesis, $h$ fulfills also the 2-term equation, we conclude from Proposition \ref{strong_unicity_of_sterling1} that $h$ is $£_1$ up to a constant, which implies that $P$ is $£_2$ up to a constant.
\end{proof}
\begin{rem}
We actually expect a slightly stronger result to be true, inasmuch as already
the duplication and 3-term relation characterize $\sterling_2$; this claim
has been verified for all primes $3<p<200$.
\end{rem}
As we can formally deduce the two equations in the proposition from the Kummer-Spence analogue and 
the Kummer-Spence analogue in turn from the Cathelineau equation $J(a,b,c)$ (because in this case the specialisation mentioned in (\ref{specialization}) is allowed),  we get
\begin{cor} The space of solutions of the Kummer-Spence analogue and the space of solutions of the Cathelineau 
equation 
are both of dimension $1$ generated by $£_2$.
\end{cor}
\begin{proof} We only need to show that if $P\in \F[T]$, assumed to be of degree less than or equal to $p-1$, setting $h=TP'$, $h$ fulfills the 2-term equation. In order to do that let $KS(a,b)$ denote the formal Kummer-Spence analogue, then taking the derivative with respect to $a$, and rewriting the equation with $h$ and finally specializing to $a=0$, we can see that, modulo the inversion formula for $h$ (which we can get directly by deriving the inversion formula for $P$), we have the identity $h(b)=h(1-b)$.
\end{proof} 



\section{Deriving functional equations : construction of the derivation map}\label{deriving}

\newcommand{\wrt}{with respect to }

The main goal of  this section is to  prove that one can pass from functional equations for polylogarithms to functional equations
for the corresponding infinitesimal polylogarithms. For this purpose we will construct a family of maps, parametrized by a given derivation, from $\B_n(F)$ to $\beta_n(F)$. The origin of such maps comes from the categorical setting which is behind the ``tangential processing'' involved in the construction of the infinitesimal polylogarithmic groups, which is to some extent discussed in \cite{bloch,cathe-sl2,
jlc-essen}, and will be treated in more detail in \cite{polyana2}.

In  subsection \ref{2}, we present the  ``derivation map'' from polylogarithmic groups to infinitesimal polylogarithmic groups. In subsection \ref{3}, we prove, as an application, that the derivation of a functional equation for any polylogarithm gives rise to a functional equation
for the corresponding infinitesimal polylogarithm, and we will show several examples. 



\subsection{From classical polylogarithmic groups to infinitesimal polylogarithmic groups}\label{2}
For the construction of the polylogarithmic groups (see section \ref{polygroups} on page \pageref{polygroups}), we gave an initial procedure for $n=2$ and an inductive procedure for higher $n$. The construction of the ``derivation map'' follows this principle.

\newsubsubsection{The case $n=2$}
\begin{lem}
Let $F$ be a field and $D \in Der_\Z(F)$ be an absolute derivation. Consider
the well-defined maps $f_D: \Z [F^{\scriptscriptstyle\bullet\bullet}] \to F[F^{\scriptscriptstyle\bullet\bullet}]$, $[a] \mapsto D(a)[a]$ and $g_D:
\Wedge^2(F^\times) \to F^\times \otimes_\Z F$, $x\wedge y \mapsto -x \otimes  \frac{D(y)}{y}+ y\otimes
\frac{D(x)}{x}$. Then the
following diagram
$$\begin{CD}
\Z [F^{\scriptscriptstyle\bullet\bullet}] @>{f_D}>> F[F^{\scriptscriptstyle\bullet\bullet}] \\
@V{\delta_2}VV @VV{\Bar{\partial}_2}V\\
\Wedge^2(F^\times) @>>{g_D}> F^\times \otimes_\Z F,
\end{CD}$$
is commutative, where $\Bar{\partial}_2([a])=\frac{1}a \otimes \frac{1}{1-a} +\frac{1}{1-a}
\otimes \frac{1}{a}$.
\end{lem}
\begin{proof} First we observe that the map $g_D$ is well defined. Indeed this is a consequence of the $d\log$ property of the map $y \mapsto \frac{D(y)}{y}$ defined on the units and of the fact that $g_D(x \otimes x)=0$ which implies that $g_D(x \wedge x)=0$. Then,  the commutativity of the diagram is a direct check. 
\end{proof}
As a direct consequence we get a map from $\ker(\delta_2)$ to $\ker(\Bar{\partial}_2)$. Similarly, 
we can obtain a map $\ker(\delta_2)$ to $\ker({\partial}_2)$ by replacing $f_D$ by $\tilde{f}_D:[a]\mapsto \frac{D(a)}{a(1-a)}[a]$ which induces a map  $\tau_{2,D}:\B_2(F)\to \beta_{2}(F)$.

\newsubsubsection{The case $n>2$} 
Suppose we have defined the ``derivation map'' $\tau_{{n-1},D}:\B_{n-1}(F)\to \beta_{n-1}(F)$ (\wrt a derivation $D$) for the level $n-1$. 
Then we can construct $\tau_{n,D}$ by induction as follows. 
\begin{prop}
Let $D\in Der_\Z(F)$ be an absolute derivation for the field $F$.
Then we have the following commutative diagram:
$$\begin{CD}
\Z [\Fdotdot] @>{\tilde{f}_D}>> F[\Fdotdot] \\
@V{\delta_n}VV @VV{{\partial}_n}V\\
\B_{n-1}(F)\otimes F^\times @>>{g_{n,D}}> \beta_{n-1}(F)\otimes_\Z F^\times \oplus \B_{n-1}(F)\otimes_\Z F
\end{CD}$$
where $\tilde{f}_D$ is defined on generators as $[a]\mapsto \frac{D(a)}{a(1-a)}[a]$, while
$g_{n,D}$ is given by 
$$g_{n,D}:\{x\}_{n-1}\otimes y \mapsto\tau_{n-1,D}\big(\{x\}_{n-1}\big)  \otimes y + \{x\}_{n-1}\otimes
\frac{D(y)}y$$
and ${\partial}_n$ by
$${\partial}_n([a])=\langle a\rangle_{n-1}\otimes a\ + \ \{a\}_{n-1}\otimes (1-a)\,.$$
\end{prop}
\begin{rem} We want to point out that despite their apparent simplicity, these crucial commutative diagrams do not show up at first sight.
\end{rem}
This induces a map from $\ker(\delta_n)$ to $\ker({\partial}_n)$ which in turn induces the desired ``derivation map''
$\tau_{n,D} :\B_n(F)\to \beta_{n}(F)$.
\begin{df} Let $F$ be a field and $D\in Der_\Z(F)$ be an absolute derivation for the field $F$. We will call the map $\tau_{n,D} :\B_n(F)\to \beta_{n}(F)$ the derivation map from $\B_n(F)$ to  $\beta_{n}(F)$, \wrt $D$. If $x$ is an element of $\B_n(F)$, the element $\tau_{n,D}(x)\in \beta_{n}(F)$ will be called the derivative of $x$ \wrt $D$.
\end{df}
As usual, if $D$ is clear from the context we will omit it. 
\begin{rem} We can notice that all the $\tau_{n,D}$, and also all the maps involved in the previous propositions, give rise to an $F$-linear map $\tau_n : Der_\Z(F) \to Hom_\Z(\B_n(F),\beta_n(F))$ for all $n\geqslant 2$.
\end{rem}
\subsection{Explicit derivation of functional equations}\label{3}
As a consequence of the previous setting we get
\begin{cor} \label{derive_for_elements}
Each element in $\ker \delta_n$ induces (many) elements in $\ker \partial_n$.
\end{cor}

The crucial main consequence is the following result.
\begin{cor} \label{derive_for_fe} Let $K$ be an arbitrary field and set $F=K(t_1,\dots,t_r)$, with $(t_1,\dots,t_r)$ a transcendence basis over $K$. Let $D\in  Der_\Z(F)$. Then any functional equation of the $n$-logarithm over $K$ induces, via the derivation map $\tau_{n,D}$, a functional equation of the infinitesimal $n$-logarithm over $K$.
\end{cor} 
\begin{proof} It is a direct consequence of the definition \ref{def_of_fe} and of the construction of  $\tau_{n,D}$.
\end{proof}
\begin{rem} Notice  that, in the above corollary, $Der_\Z(F)\neq 0$ since $Der_K(F)\neq 0$, at least if $r\geqslant 1$. In practice it could be interesting to have a differential basis, and thus we can assume that if $K$ is of characteristic $p$ then $(t_1,\dots,t_r)$ is a $p$-basis over $K$.
\end{rem}

It is a priori not clear that the procedure gives non-trivial equations, but the following examples show that it is actually the case:

\begin{ex} The first example is taken from \cite{cathe-pol} and it retrieves the 4-term relation from
the 5-term relation (\ref{five2}) by applying the above procedure with
$$D=a(1-a)\frac{\partial}{\partial a} +b(1-b) \frac{\partial}{\partial b}\,,$$
assuming that $F=K(a,b)$ with $a,b$ indeterminates over the field $K$, and that $\frac{\partial}{\partial a}$ and $\frac{\partial}{\partial b}$ are the usual partial derivatives.
\end{ex}

The following proposition gives a partial answer to Cathelineau's question
concerning the relationship of his 22-term equation for $\,d\pol_3\,$ and
Goncharov's equation (\ref{goncha}) for $\,\pol_3\,$ (with the same number of terms).
It is a consequence of the previous results but can also be verified directly.
\begin{prop} \begin{enumerate}
\item The infinitesimal functional equation below, which is derived from the Goncharov 
functional equation for the trilogarithm
is zero in  $\,\beta_3(F)\,$. 
\item If  $\,F \supset \Q$, the Goncharov equation is expressible in terms of  an $F$-linear combination of  $\,J(a,b,c)\,$.
\end{enumerate}
\end{prop}
We give an example of such a derived version in the case $F=K(a,b,c)$ with $a,b,c$ indeterminates over the field $K$, applying the above procedure with
$$D=a(1-a)\frac{\partial}{\partial a} +b(1-b) \frac{\partial}{\partial b}+c(1-c) \frac{\partial}{\partial c}\,$$
to the equation stated in (\ref{goncha}). 
Let us set
\begin{align*} 
\varphi(a,b,c)=&[a]-\frac {( b- 1)( a- 1 )}{ab-1}\left( \left[-  \frac {b( a- 1 )}{b - 1}\right]
+\left[ - \frac {a( b- 1
)}{ a- 1 }\right]\right)\\
&+ \frac {( c^{2}b + c
b^{2} -3cb+ 1)}{cb-1} \left[\frac { a- 1}{a
bc - 1}\right]-\frac {(abc - a - b - c + 2)}{cb-1}\left[ \frac {cb(a - 1)}{ab
c - 1}\right]\\
& - 
 \frac {(a +b- 2)}{ab-1}[ab]-\frac {(a^{2}bc - 2abc + 
b +c- 1 )( a- 1)}{(ac-1)(ab-1)} \left[- 
 \frac {a( c- 1)(
b - 1)}{( a- 1)(abc - 1)}\right]\,.
\end{align*}
Then, modulo the inversion formula,
$$\varphi(a,b,c)+\varphi(b,c,a)+\varphi(c,a,b)-\frac {(a + b + c - 3)}{abc - 1}[abc]$$
is the differential of the Goncharov equation and vanishes in
$\beta_3(F)$ by virtue of Corollary \ref{derive_for_elements}.\\[0.1cm] 

\begin{obs} We should notice that we have not yet  proved that the infinitesimal Goncharov equation also holds in characteristic $p$ and to know that this equation is expressible in terms of  an $F$-linear combination of  $\,J(a,b,c)\,$ is not enough to ensure this (unless we know that this linear combination is independant of $F$). It will be seen in the next section that it is the case, at least if we see $\sterling_2$ as a function from $\Z/p$ to $\Z/p$.
\end{obs}




\section{Reduction of functional equations mod $p$ via the $p$-adic realm}\label{red_mod_p}


\newcommand{\fe}{{functional equation}}
\newcommand{\pl}{{polylogarithm}}
\newcommand{\loga}{{logarithm}}
\newcommand{\polylog}{{\mathcal D}}
\newcommand{\field}{{\Q_p(t_1,\dots,t_r)}}
\newcommand{\Der}{{\mbox{Der}}}
\newcommand{\ratfctfield}{{\Q(t_1,\dots,t_r)}} 

In this section, we want to prove the following statements (which are made more precise below):
\begin{enum}
\item \label{1)}Each \fe\ for the classical $n$-logarithm $\polylog_n$ induces a functional equation
for certain $p$-adic $n$-logarithm functions (those which satisfy Wojtkowiak's $p$-adic version 
of Zagier's criterion).\\
\item \label{2)} Each \fe\ for the classical $n$-logarithm induces a functional equation for the 
infinitesimal $n$-logarithm (via the derivation procedure given in the previous section). A similar statement holds 
for the $p$-adic case.\\
\item \label{3)}Each \fe\ for the infinitesimal $n$-logarithm $d\polylog_n$ induces a \fe\ for the
corresponding $p$-adic infinitesimal $n$-\loga\ denoted $DF_n$ (see Definition  \ref{besserdef}).\\
\item \label{4)}Each ``good $\Q_p$-specialization'', as defined in (\ref{good_reduction}) below, of a \fe\ for the $p$-adic infinitesimal \pl\ induces a \fe\ (in  the weak sense) for the 
finite $(n-1)$-\loga.
\end{enum}

Combining the four statements, we arrive at the somewhat more surprising statement:

\medskip\noindent
{\bf Surprise:} {\it Each \fe\ for the classical $n$-\loga\ induces a \fe\ for the finite $(n-1)$-\loga.}\\

Throughout this section, we denote by $Li_n(z)$ Coleman's $p$-adic $n$-logarithm \cite{coleman}.
Let us first look for the $p$-adic combinations which should play the same
role as the modified polylogarithms $\polylog_n$.

\begin{rem} The inversion relation \big(in its clean form $P_n(z)=(-1)^{n-1}P_n(1/z)$\big) for a combination $P_n(z)=\sum_{k=0}^{n-1} a_k \log^k(z) Li_{n-k}(z)$
is equivalent to the following condition on the coefficients:
\begin{equation}\label{inv-cond}
\sum_{k=0}^{n-1} \frac{a_k}{(n-k)!} = 0
\end{equation}
(cf. \cite{wojt2}, Lemma 4.2). Since the inversion
relation is in the kernel of $\partial_n$, we can restrict our investigations to combinations $P_n(z)$ 
satisfying those conditions.
\end{rem}
While one needs to work harder in the ``classical'' case to find functions 
which satisfy cleanly their functional equations, it turns out that in the 
$p$-adic case the above condition is already good enough, and
we can state the above claim \ref{1)} more precisely as

\begin{prop} {\rm (Wojtkowiak, \cite{wojt2}, Proposition 4.4)}\\
Let $\xi\in \ker\delta_{n,\field}$. Then each admissible $\C_p$-specialization of $\xi$  is mapped to a constant by the $p$-adic functions
\begin{equation} \label{polycomb}
P_n(z)=\sum_{k=0}^{n-1} a_k \log^k(z) Li_{n-k}(z)\,,\end{equation}
if the coefficients satisfy condition  {\rm (\ref{inv-cond})}.
\end{prop}

This motivates the following definition:
\begin{df}
A linear combination of $p$-adic polylogarithms of the form {\rm (\ref{polycomb})}
whose coefficients satisfy {\rm (\ref{inv-cond})} is called a {\rm
clean $p$-adic polylogarithm}.
\end{df}

\begin{rem}
\begin{enumerate}
\item For $n=2$, there is, up to a multiplicative
constant, only one clean $p$-adic 2-logarithm $P_2$ satisfying (\ref{inv-cond}).
\item The original statement was actually somewhat stronger:  $\field$ was replaced
by $\C_p(t)$, where
$\C_p$ denotes a completion of an algebraic closure of $\Q_p$.
\end{enumerate}
\end{rem}

\noindent The claim in \ref{2)} follows immediately from the ``derivation map'' in \S\ref{deriving}.\\
\

Before we show a more precise version of \ref{3)}  by imitating Proposition 7 of \cite{cathe-pol}, we 
state our intermediate goal: We are looking for a morphism $F[F]\to F_0$, where $F=\C_p(z)$ and $F_0=\C_p$.
More precisely, we want to have a family of morphisms $(DP_n)_{n\geqslant 2}$ on $\beta_n(\C_p)$ expressed in terms 
of the differential operator $D=z(1-z)\frac d {dz} $ and some clean $p$-adic polylogarithms $P_n$. There are many 
candidates:

\begin{prop}\label{vanishing_for_p_adic_diff}
Let $(P_n)_{n\geqslant 2}$ be a family of clean $p$-adic polylogarithms 
such that for $n\geqslant 3$
\begin{equation} \label{diff_eq}
DP_n(z)=\lambda_n\, (1-z) P_{n-1}(z) \ +\ \mu_n\, \log(z)\, DP_{n-1}(z)\,,
\end{equation}
for some $\lambda_n$, $\mu_n\in \C_p^\times$. \\
Then, for any $n$, $DP_n$ defines a morphism on $\beta_n(\C_p)$.
\end{prop}
\begin{proof}
$P_n$ is defined on $\B_n(\C_p)$ by assumption. For $n=2$, we have seen that the function is essentially unique:
$$P_2(z)=-2\,Li_2(z)+ \log(z) Li_1(z)\,,$$
and the resulting infinitesimal dilogarithm
$$DP_2(z)=(1-z)\log(1-z) + z\log(z)$$
vanishes on $r_2(\C_p)$ (due to Proposition 2.8, it is enough to check that it vanishes on the four term relation,
which is straightforward).\\
Now suppose the claim is true for $n-1$. 
The maps $\,DP_{n-1}\otimes \log:\beta_{n-1}(\C_p)\otimes \C_p^\times\to \C_p$, $x\langle y\rangle_{n-1}\otimes z\mapsto xDP_{n-1}(y)\log(z)$ resp. $\,P_{n-1}\otimes Id:\B_{n-1}(\C_p)\otimes \C_p\to \C_p$, $\{y\}_{n-1}\otimes z\mapsto zP_{n-1}(y)$, are well-defined by the inductive assumption resp. by assumption ($P_{n-1}$ is clean).
Furthermore, an element $\xi\in r_n(\C_p)$ lies in the kernel of each of the ``components'' of $\partial_n$, say
$\partial_n':\C_p[\C_p]\to \beta_{n-1}(\C_p)\otimes \C_p^\times$ and
$\partial_n'':\C_p[\C_p]\to \B_{n-1}(\C_p)\otimes \C_p$, and therefore
$$\Big(\mu_n\, DP_{n-1}\otimes\log\ +\ \lambda_n\ P_{n-1}\otimes Id\Big)\big(\partial_n\xi) =\Big(\mu_n\, DP_{n-1}\otimes\log\circ\partial_n' +\lambda_n\ P_{n-1}\otimes Id\circ \partial_n''\Big)\xi=0\,,$$
which shows that the function defined by $(\ref{diff_eq})$ can be linearly extended to a well-defined function on $\beta_n(\C_p)$.
\end{proof}

\begin{df}\label{besserdef} 
Besser's $p$-adic $n$-logarithm is defined  as
\begin{equation}\label{besser-family} 
F_n(z)=\sum_{k=0}^{n-1} a_{k,n} \log^k(z) Li_{n-k}(z)
\end{equation}
with
$$ a_{k,n}=\frac{(-1)^k}{k!} (k-n)\,.$$
We will call $DF_n$ the {\rm distinguished infinitesimal $p$-adic $n$-logarithm}.
\end{df}

\begin{prop}[Existence]
There exist families of clean $p$-adic polylogarithms satisfying {\rm (\ref{diff_eq})} for some $\lambda_n$, $\mu_n\in \C_p^\times$.\\
In particular, Besser's family {\rm (\ref{besser-family})}
satisfies {\rm (\ref{diff_eq})} with $(\lambda_n,\mu_n)=(\frac 1{n-1},-\frac 1{n-1})$, $n\geqslant 3$.\\
There are many other possibilities.
\end{prop}
\begin{proof}
Again, the case $n=2$ gives the unique choice for $P_2$ (up to multiplicative constant).\\
Inductively, starting from $P_{n-1}$ and $DP_{n-1}$, one can form an arbitrary linear combination of them using
$\lambda_n$ and $\mu_n$ which gives a candidate for $DP_n$, with coefficients $b_{k,n}$, say; a subsequent 
``integration'' (putting $a_{0,n}=-n$ and successively $a_{k+1,n}=-n(b_{kn}-a_{kn})/(k+1)$, $k=0,\dots,n-2$) provides a candidate 
$P_n$ whose coefficients $a_{kn}$ have to satisfy the further condition (\ref{inv-cond})---this gives a linear restriction on the possible $(\lambda_n,\mu_n)$ at each step. We thus obtain inductively
an extra degree of freedom at each level.\\
For example, normalizing $P_n(z)$ such that  $a_0-n$, we obtain successively
$$\lambda_3-\mu_3=1\,,\qquad \lambda_4-\mu_4=\frac 1{2-\lambda_3}\ ,\ {\rm etc.}$$
It remains to check that Besser's choice (\ref{besser-family}) does satisfy
\begin{equation} \label{diff-besser}
(n-1)\,DF_n(z)= (1-z)\,F_{n-1}(z) - \log(z) DF_{n-1}(z)
\end{equation}
which is straightforward. Also, the $a_{k,n}$ satisfy condition (\ref{inv-cond}) since
$$-\sum_{k=0}^{n-1}\frac{(-1)^k}{k!(n-k)!} (n-k)= \frac 1{(n-1)!} (1-1)^{n-1}=0\,.$$
\end{proof}
\begin{rem} 
\begin{enumerate}
\item Writing $\Phi_n(z)= (n-1)! F_n(z)$ and noticing that $(1-z)= D\log(z)$, we can reformulate
(\ref{diff-besser}) more suggestively, using the ad-hoc convention $D^-(a\otimes b):=D(a)b - a D(b)$, as
$$D\Phi_n(z) = D^-\big(\log(z)\otimes \Phi_{n-1}(z)\big)\,.$$
\item
We have just seen that, a priori, there are many choices for the $P_n$ individually,
but the condition that the morphisms at level $n$ and $n-1$ be linked via the condition
$\rho\,DP_n(z)=(1-z) P_{n-1}(z)-\log(z) DP_{n-1}(z)$ for some $\rho\in\C_p$ provides us with a unique function, up to a 
multiplicative factor,  the condition  (\ref{inv-cond}) still being true for $P_n$. We have not found a ``natural''
justification for the condition (\ref{diff-besser}), though. A normalization condition for the above $P_n$ is then $a_{0,n}+a_{1,n}=-1$ which entails
$\rho =n-1$. 
The resulting family coincides with Besser's functions (\ref{besser-family})---his choice of coefficients was forced by two
rather natural requirements: first, a certain $p$-adic power series expansion becomes independent of the ``direction'' in which to expand; second,
one retrieves the finite $(n-1)$-\loga\ by
reducing $DF_n$ mod $p^n$ (or, more precisely, reducing $p^{1-n}DF_n$ mod $p$) on 
elements in $\Z_p^\times\cap(1-\Z_p)^\times\subset \C_p$ (for an improved statement of this and of the following
theorem cf.~\cite{Besser}). 
\end{enumerate}
\end{rem}

The $F_n$ can be characterized by the following

\begin{thm}\label{besserthm} {\rm (Besser, \cite{Besser}, Theorem 1.1)}\\
Let $X=\{z\in \Z_p \, , \,|z|=|1-z|=1\}$. For $p>n+1$, one has  $DF_n(\Z_p)\subset p^{n-1}\Z_p$, 
and for $z\in X$: 
\[p^{1-n}DF_n(z)\equiv \pounds_{n-1}(z) \pmod p\, .\]
The choice of coefficients (in $\Q$) for $F_n$ is unique for a clean $p$-adic polylogarithm which satisfies the above property for all $p>n+1$.
\end{thm}

In order to formulate the subsequent statements conveniently, we introduce the following notion:
\begin{df}\label{good_reduction}
A {\rm good $\Q_p$-specialization for} $\sum n_i[x_i]\in F[F]$, $F\subset\field $, is a 
family of numbers $u_j\in \Q_p$, $j=1,\dots,r$, such that the images of $n_i=n_i(t_1,\dots,t_r)$, $x_i=x_i(t_1,\dots,t_r)$ and $1-x_i=1-x_i(t_1,\dots,t_r)$ under the specialization map $t_j\mapsto u_j$, $j=1,\dots,r$, are in $\Z_p^\times$.
\end{df}

The virtue of a good $\Q_p$-specialization lies in the fact that we can reduce it modulo $p\Z_p$. As we can notice, a good $\Q_p$-specialization is, in particular, an admissible $\Q_p$-specialization.
Now, putting Proposition \ref{vanishing_for_p_adic_diff} and Theorem \ref{besserthm} together, we can make \ref{4)} more precise:
\begin{cor}
Let $n\geqslant 2$, $p>n+1$, and $\eta\in \ker \partial_{n,\field}$. Then we have
\begin{itemize}
\item[a)] For each admissible $\C_p$-specialization $\eta^{\rm spec}$ for $\eta$, $DF_n(\eta^{\rm spec})=0$.
\item[b)] For each good $\Q_p$-specialization $\eta^{\rm spec}$ for $\eta$, the reduction
mod $p$ gives
$$\pounds_{n-1}(\eta^{\rm spec})\equiv 0 \pmod p\,.$$
\end{itemize}
\end{cor}

\begin{proof}
The infinitesimal \pl\ $DF_n$ vanishes on $\eta$ by Proposition \ref{vanishing_for_p_adic_diff}, and reducing mod $p$ obviously
conserves this vanishing property. Besser's result now says that the reduction of
$p^{1-n}DF_n(\eta^{\rm spec})$ is equal to $\pounds_{n-1}\big(\eta^{\rm spec} \pmod p\big)$.
\end{proof}

Going even one step further, we can state a more precise version of the above ``surprise'':
\begin{cor}\label{surprise}
Let $n\geqslant 2$, $p>n+1$, and $\xi\in \ker \delta_{n,\ratfctfield}$. Then we have
\begin{itemize}
\item[a)] For each admissible $\C$-specialization resp. $\C_p$-specialization $\xi^{\rm spec}$ for $\xi$, the quantities  $\polylog_n(\xi^{\rm spec})$ resp.  $F_n(\xi^{\rm spec})$ are constants.
\item[b)] For each absolute derivation $\Delta\in \Der_\Z(\ratfctfield)$, $\xi$ induces
$\xi_\Delta\in \ker \partial_{n,\ratfctfield}$, and therefore, for each admissible $\C$-specialization resp. 
 $\C_p$-specialization,
$$d\polylog_n(\xi_\Delta)=0,\quad \mbox{resp.} \quad DF_n(\xi_\Delta)=0\,.$$
\item[c)] For each good $\Q_p$-specialization $\xi_\Delta^{\rm spec}$ for $\xi_\Delta$,
the reduction mod $p$ gives
$$\pounds_{n-1}(\xi_\Delta^{\rm spec})\equiv 0 \pmod p\,.$$
\end{itemize}
\end{cor}

\begin{proof} \begin{enumerate} \item[a)]
 Follows from Zagier \cite{z1}  and Wojtkowiak \cite{wojt2}, respectively.  
\item[b)] This follows via the ``derivation map'' (see §\ref{deriving}). 
\item[c)] $0=p^{1-n}DF_n(\xi_\Delta)\equiv \pounds_{n-1}(\xi_\Delta^{\rm spec})$.\end{enumerate}
\end{proof}

Alas, although being quite powerful, the above strategy does not give the full answer to our problem.

\begin{rem}
\begin{enumerate}
\item The virtues of the procedure described above lie in its generality: we do not need
to (find and) prove functional equations for ($p$-adic) infinitesimal or finite polylogs,
since they ``drop out'' using the machinery.\\
\item The drawbacks of the machinery lie in its lack of control: 
\begin{enumerate}
\item We do not get the \fe s as polynomial identities but only ``on points'', i.e.~in the form of
(good) specializations.\\
\item A more mundane reason for proving functional equations for $\pounds_n$ in
the strong sense is
the fact that all the ones which have occurred in our investigations are not only true for
$\F_p$ but actually hold more generally for any field of characteristic $p$.\\
\item (a minor point, given the range in which we mostly work) We need to assume that $p>n+1$.\\
This restriction is not (always) necessary for the polynomial identities to hold: there are examples
of equations for $\sterling_3$ which are still true in characteristic 3.
\end{enumerate}
\end{enumerate}
\end{rem}

In summary, there are still plenty of reasons which leave us with the task of finding proofs of  
\fe s for the finite polylogarithms. The final section will therefore be dedicated to this issue.


\part{The Main Proofs}

\section{Proofs of functional equations over fields of characteristic $p$.}\label{theproofs}
\subsection{Straightforward demonstrations}
\begin{proof} (of Proposition \ref{generic})\\
\begin{enumerate}
\item The inversion relation can be checked via a straightforward algebraic manipulation.
\item In order to prove the distribution relation, let us fix a primitive $m$th root of unity $\,\zeta\,$. 
Dividing both sides by $\,m^n\,$ and developing the fraction into
a (finite) series leaves us to prove:
\begin{align*}
\sum_{k=1}^{p-1}\frac{T^{km}}{(km)^n} &=\frac 1 m\sum_{\zeta^m=1}
\Big(1+(\zeta T)^p+(\zeta T)^{2p}+\dots+(\zeta T)^{(m-1)p}\Big)\sum_{k=1}^{p-1}
\frac{(\zeta T)^k}{k^n}\\
&=\frac 1 m\sum_{k=1}^{p-1}\frac{1}{k^n}\sum_{\zeta^m=1}\Big((\zeta T)^k+(\zeta T)^{p+k}+(\zeta T)^{2p+k}+\dots+(\zeta T)^{(m-1)p+k}\Big)\\
&=\frac 1 m\sum_{k=1}^{p-1}\sum_{\zeta^m=1}\bigg(\frac{(\zeta T)^k}{k^n}+\frac{(\zeta T)^{p+k}}{(p+k)^n}+\frac{(\zeta T)^{2p+k}}{(2p+k)^n}+\dots+\frac{(\zeta T)^{(m-1)p+k}}{\big((m-1)p+k\big)^n}\bigg)\\
&=\frac 1 m\sum_{\textstyle{r=1\atop p\hskip -2pt\not\hskip 2pt|r}}^{pm-1}\Big(\sum_{\zeta^m=1}\zeta^r\Big)\frac{T^r}{r^n}\,,
\end{align*}
and this is true due to the character relations
$$\sum_{\zeta^m=1}\zeta^r=\begin{cases} m,&\text{if }m|r\\ 0,& \text{otherwise}.\end{cases}
$$

\item (Proof of the special values) 
$\wt{{£}_{n}}(1)=0\,$ if $\,(p-1)\hskip -4pt\not|n\,$ follows from the well-known fact that 
$\sum_{k=0}^{p-1} P(k)=0$ for any polynomial $P\in\Z/p\Z[x]$ of degree $\leqslant p-2$
(here we apply it to the monomials $x,\dots,x^{p-2}$), the statement for $\,(p-1)|n$ being
obvious. \\
The assertion for $\wt{{£}_{2n}}(-1)=0$ is a direct consequence of the inversion relation.

To prove the last formula of Proposition \ref{generic} we only need to take $m=2n$ (the odd values correspond to the above identities). For this, 
one can use the special case $a=2$, in \cite{rib}(Proposition (5B), p.108),
$p-1\not| 2n$:
\begin{equation}\label{ribcong}
(1-2^{2n})B_{2n}\equiv 2n\,2^{2n-1} \sum_{1\leqslant j<\frac p 2}\frac 1{j^{1-2n}}
\pmod p 
\end{equation}
and the fact that $\sterling_{1-2n}(-1)$ is equal to the sum in (\ref{ribcong}): rewriting 
$$\wt{{£}_{p-2n}}(-1)=\wt{{£}_{1-2n}}(-1)=\sum_{j=1}^{(p-1)/2}(2j)^{2n-1}-\sum_{j=1}^{(p-1)/2}(2j-1)^{2n-1}
=2\sum_{j=1}^{(p-1)/2}(2j)^{2n-1}-\sum_{j=1}^{p-1} j^{2n-1}\,,$$
one sees that the first sum is equal to $2^{2n}$ times the sum in (\ref{ribcong}), while the
second one equals $-\wt{{£}_{1-2n}}(1)$ and therefore is zero (for $0<n<\frac{p-1}2$)  by the above special value.
\end{enumerate}
\end{proof}

\subsection{A recipe for proving functional equations}
Let $R$ be a domain of characteristic $p$. In order to show that a polynomial $Q(T)\in  R[T]$ is 
zero, we divide it into three parts:
$$Q(T)=Q(0) +Q_1(T)+Q_2(T^p)\,,$$
where $Q_1(T)$ involves only powers of $T$ which are not divisible by $p$. Then we verify separately that
$Q_2(T^p)$ and the constant $Q(0)$ vanish and that $\frac d {dT} Q_1(T)$ is zero as well.
We can iterate this procedure in an obvious way.
\\

\begin{proof}[Proof of Proposition \ref{L1.1}]
\begin{enumerate}
\item 
We will apply the recipe above. We have
\begin{align*}
\frac{d}{dT}£_{1}(1-T)&=-\frac{1}{1-T}£_{0}(1-T),\\
&=\frac{1}{T}£_{0}(T) \quad {\rm by }\; (\ref{L0}),\\
&=\frac{d}{dT}£_{1}(T),
\end{align*}
and as the degree of either polynomials is less than $\,p-1\,$, we conclude that
$\,£_{1}(T)=£_{1}(1-T)+c\,$ where $\,c\,$ is a constant. \\
This, in turn, implies that $\,2c=0\,$ (specialize $\,T=0\,$ and $\,T=1\,$, 
respectively), and therefore we get as a by-product $\,£_{1}(1)=£_{1}(0)=0\,$
(in characteristic $\not= 2$).\\

\item 
The following proof is a slight variation of the recipe, in that it uses two 
iterated derivatives.\\
Denote by $\,\partial_x\,$ and $\,\partial_y\,$ the derivatives to respect to $\,x\,$ and $\,y\,$. 
We can check, using the differential equation for $\,\sterling_1\,$ and the
rational expression (..) for $\,\sterling_0\,$, that
$$\partial_y \partial_x H(x,y,s)=\frac{1-y^p-x^p+s^p}{(1-y-x+s)^2},$$
which is an expression that is symmetric in $\,x\,$ and $\,y\,$. Thus
$$\partial_y \partial_x (H(x,y,s)-H(y,x,s))=0.$$
But the maximum degree for each indeterminate in the polynomial $\,H(x,y,s)\,$ 
is never greater than $\,p-1\,$, and as a consequence the above identity implies that
$$H(x,y,s)-H(y,x,s)=R_0(s)+R_1(s) x+R_2(s)y,$$
where $\,R_0,R_1,R_2 \in \F[s]\,$.
But setting $\,x=y\,$ implies both
$\,R_0=0\,$ and $R_1+R_2=0$, and the construction of $\,R_1\,$ and $\,R_2\,$ shows directly that they are
both zero \big(the coefficients of $x$ and $y$ in $H(x,y,s)$ are both equal to $\sum_{k=0}^{p-2}
(-s)^k$\big).
\end{enumerate}
\end{proof}

\begin{proof}[Proof of Proposition \ref{L2.0}]
\begin{enumerate}
\item 
Set
$$E(T)=£_{2}(1-T)-£_{2}(T)+T^p £_{2}(1-\frac{1}{T}).$$
We want to prove that $\,E\,$ is 0 in $\,\F[T]\,$. 
Computing $\,\frac{d}{dT}E\,$ we get
$$\frac{d}{dT}E(T)=-\frac{1}{1-T}£_1(1-T)-\frac{1}{T}£_1(T)+\frac{T^{p-1}}{T-1}£_1(1-\frac{1}{T}).$$
But by Proposition \ref{L1.1}, $\,£_1(1-T)=£_1(T)\,$ and $\,£_1(1-\frac{1}{T})=£_1(\frac{1}{T})\,$. Moreover by the inversion formula (see Proposition \ref{generic}) we have $\,£_1(\frac{1}{T})=-\frac{1}{T^p}£_1(T)\,$. Hence,
\begin{align*}
\frac{d}{dT}E(T)&=-\frac{1}{1-T}£_1(T)-\frac{1}{T}£_1(T)-\frac{1}{(T-1)T}£_1(T),\\
&=0.
\end{align*}
As $\,E(0)=0\,$ and $\,\deg(E)\leqslant p\,$, we know that $\,E(T)=cT^p\,$ and therefore $T^p E(\frac1T)=c$,
but using the inversion relation one sees that $T^p E(\frac1T)=E(T)$, which implies $c=0$.
\end{enumerate}
\end{proof}
\begin{rem}
A different way to prove that $\,c=0\,$: For this we look directly at  $\,E(T)\,$ and try to compute this coefficient which can only appear in the expression
$$T^p £_{2}(1-\frac{1}{T})=\sum_{i=1}^{p-1}\frac{T^{p-i}(T-1)^i}{i^2}.$$
But, for each $\,i\,$, the coefficient of $\,T^p\,$ is $\,\frac{1}{i^2}\,$, and thus $\,c=£_{2}(1)=0\,$.\\
\end{rem}

\begin{proof}[Proof of Theorem \ref{L2.1}]
The strategy of proof could be summarized as follows:\\
\begin{itemize}
\item[(i)] Prove that $\partial_c \widehat{£_2}(J(a,b,c))=0$ in $\F[a,b,c]$.
\item[(ii)] Prove that $\widehat{£_2}(J(a,b,0))=0$ in $\F[a,b]$.
\item[(iii)] Prove that the coefficient of $c^p$ in $\widehat{£_2}(J(a,b,c))$ is 0.
\end{itemize}
For the proof of this functional equation we will need several preliminary formulas. First we will use these two relations, in $\F[x,y]$, coming from the 4-term equation for $£_1$
\begin{align}
(1-y)^p £_1 \left( \frac{x}{1-y}\right) &= £_1(x)+(1-x)^p £_1\left( \frac{y}{1-x}\right)-£_1(y), \label{R1} \\
£_1(y) -£_1(x)&= (1-x)^p £_1 \left( \frac{1-y}{1-x}\right) +x^p £_1\left( \frac{y}{x}\right).\label{R2}
\end{align}
We use implicitly the following formal derivation rules, where $t$ is an indeterminate and $\lambda$ a constant independent of $t$:
\begin{align*} 
\frac{d}{dt}£_2\big(\lambda (1-t)\big) &= -\frac{1}{1-t}\, £_1(\lambda (1-t)),\\
\frac{d}{dt}£_2\left(\lambda \Big( 1-\frac{1}{t}\Big)\right) &= -\frac{1}{t(1-t)} \,£_1\left(\lambda \Big( 1-\frac{1}{t}\Big)\right).
\end{align*}
We also point out that the following simple formula will be often used:
$$\frac{1}{t}+\frac{1}{1-t}=\frac{1}{t(1-t)}.$$
For the convenience of the reader we will give detailled computations in order to make checking almost straightforward.\\

Let's first split $\widehat{£_2}(J(a,b,c))$ into six pieces to facilitate the identification of the cancellation in the forthcoming computations:
\begin{align*}
A_1 &= c^p£_2(a) - 
c^p £_2(b) +(a-b+1)^p£_2( c)\\[0.1cm]
&+(1-c)^p£_2( 1-a) -(1-c)^p£_2( 1-b) 
+ (b-a)^p£_2( 1-c),\\[0.1cm]
A_2 &=-a^p£_2\left( \frac{c}{a}\right)  
+b^p£_2\left( \frac{c}{b}\right) 
+c^pa^p£_2\left( \frac{b}{a}\right)\\[0.1cm]
&-(1-a)^p£_2\left( \frac{1-c}{1-a}\right)
+(1-b)^p£_2\left( \frac{1-c}{1-b}\right) 
+c^p(1-a)^p£_2\left( \frac{1-b}{1-a}\right),\\[0.1cm]
A_3 &=c^p(1-a)^p£_2\left( \frac{a(1-c)}{c(1-a)}\right) 
-c^p(1-b)^p£_2\left( \frac{b(1-c)}{c(1-b)}\right),\\[0.1cm]
A_4 &=-b^p£_2\left( \frac{ca}{b}\right) 
-(1-b)^p£_2\left( \frac{c(1-a)}{1-b}\right), \\[0.1cm]
A_5 &=-(a-b)^p£_2\left( \frac{(1-c)a}{a-b}\right) 
-(b-a)^p£_2\left( \frac{(1-c)(1-a)}{b-a}\right) \\[0.1cm]
&+c^p(a-b)^p£_2\left( \frac{(1-c)b}{c(a-b)}\right) 
+c^p(b-a)^p£_2\left( \frac{(1-c)(1-b)}{c(b-a)}\right),\\
A_6 &=(1-c)^pa^p£_2\left( \frac{a-b}{a}\right) 
+(1-c)^p(1-a)^p£_2\left( \frac{b-a}{1-a}\right).
\end{align*}
Set $d=\frac{\partial}{\partial c}$.\\
{\it First step}: {\sl prove that $\sum dA_i=0$}.\\
It is immediate that $dA_6=0$. Using the rules described above, we get the following equalities:
\begin{align*}
dA_1 &= \frac{1}{c} £_1(c) + \frac{(a-b)^p}{c(1-c)}£_1(c),\\
dA_2 &= -\frac{a^p}{c} £_1\left( \frac{c}{a}\right)+\frac{b^p}{c} £_1 \left( \frac{c}{b} \right)\\
&+\frac{(1-a)^p}{1-c}£_1\left(\frac{1-c}{1-a}\right)-\frac{(1-b)^p}{1-c} £_1\left(\frac{1-c}{1-b}\right),\\
dA_3 &= -\frac{c^p(1-a)^p}{c(1-c)} £_1\left( \frac{a(1-c)}{c(1-a)}\right)+\frac{c^p (1-b)^p}{c(1-c)} £_1 \left( \frac{b(1-c)}{c(1-b)}\right),\\
dA_4 &= -\frac{b^p}{c} £_1 \left( \frac{ca}{b}\right)-\frac{(1-b)^p}{c} £_1\left( \frac{c(1-a)}{1-b}\right),\\
dA_5 &= \frac{(a-b)^p}{1-c} £_1\left( \frac{(1-c)a}{a-b}\right)+\frac{(b-a)^p}{1-c}£_1\left( \frac{(1-c)(1-a)}{b-a}\right)\\
&-\frac{c^p(a-b)^p}{c(1-c)}£_1\left( \frac{(1-c)b}{c(a-b)}\right) -\frac{c^p(b-a)^p}{c(1-c)} £_1 \left( \frac{(1-c)(1-b)}{c(b-a)}\right).
\end{align*}
Then, applying consecutively (\ref{R1}) to $dA_5$, with
 $x=1-c$, $y=\frac{b}{a}$, with $x=1-\frac{1}{c}$, $y= \frac{b}{a}$, and with 
$x=1-\frac{1}{c}$, $y= \frac{1-a}{1-b}$, and to $dA_3$ with $x=1-\frac{1}{c}$,
 $y=\frac{1}{a}$, and using (\ref{R2}) for simplification as well as the basic 
relations for $£_1$, we get
\begin{align*}
dA_4+dA_5+dA_3 &= -\frac{£_1(b)}{1-c}+\frac{£_1(a)}{1-c}\\
&+\frac{c^p}{c(1-c)}\left(£_1\left(\frac{b}{c}\right)-£_1\left(\frac{a}{c}\right)\right)-£_1(c)+\frac{(b-a)^p}{c(1-c)}£_1(c),
\end{align*}
then
$$dA_1+dA_4+dA_5+dA_3 =-\frac{£_1(b)}{1-c}+\frac{£_1(a)}{1-c}+\frac{c^p}{c(1-c)}\left(£_1\left(\frac{b}{c}\right)-£_1\left(\frac{a}{c}\right)\right).$$
It remains to transform $dA_2$, but  using (\ref{R2}), we have e.g.
\begin{align*}
(1-b)^p £_1\left(\frac{1-c}{1-b}\right) &=£_1(c)-£_1(b)-b^p£_1\left(\frac{c}{b}\right),\\
\end{align*}
then
$$dA_2=\frac{£_1(b)}{1-c}-\frac{£_1(a)}{1-c}+\frac{b^p}{c(1-c)}£_1\left(\frac{c}{b}\right)-\frac{a^p}{c(1-c)}£_1\left(\frac{c}{a}\right).$$
Now by  invoking the inversion formula we see that
$$\sum_{i=1}^5 dA_i =0.$$

\noindent{\it Second step}: {\sl Prove that the relation is true for $c=0$}.\\
Putting $c=0$ in $\sum_{i=1}^6 A_i$ gives 
\begin{align*}
£_2(1-a)-£_2(1-b)-(1-a)^p£_2\left(\frac{1}{1-a}\right)+(1-b)^p £_2\left(\frac{1}{1-b}\right)\\
a^p£_2\left(\frac{a-b}{a}\right)+(1-a)^p£_2\left(\frac{b-a}{1-a}\right)-(a-b)^p£_2\left(\frac{a}{a-b}\right)-(b-a)^p£_2\left(\frac{1-a}{b-a}\right)
\end{align*}
and applying the inversion formula for $£_2$ we get 0.\\

\noindent{\it Third step}: {\sl Prove that the coefficient of $c^p$ is 0}.\\
Notice first that if $\lambda$ is an expression independent of $c$, then the coefficient of $c^p$ in the sum $\sum_{i=1}^{p-1}\frac{\lambda^i}{i^2}c^{p-i}(1-c)^i$ is $£_2(-\lambda)$. Using this fact, we can see that the coefficient of $c^p$ in the expression $\sum_{i=1}^6 A_i$ is given by
\begin{align*}
&£_2(a)-£_2(1-a)+(1-a)^p£_2\left(\frac{-a}{1-a}\right)\\
&-£_2(b)+£_2(1-b)-(1-b)^p£_2\left(\frac{-b}{1-b}\right)\\
&+a^p£_2\left(\frac{a}{b}\right)-a^p£_2\left(\frac{a-b}{a}\right)+(a-b)^p£_2\left(\frac{-b}{a-b}\right)\\
&+(1-a)^p£_2\left(\frac{1-b}{1-a}\right)-(1-a)^p £_2 \left(\frac{b-a}{1-a}\right)+(b-a)^p£_2\left(-\frac{1-b}{b-a}\right).
\end{align*}
But each of the previous lines are 0 by using the 3-term functional equation of $£_2$ (see Proposition \ref{L2.1} 1. ) and this completes the proof of the 22-term functional equation for $£_2$.
\end{proof}

\begin{rem} 
We want to stress some more structural properties in the rather computational parts of the 
previous proof---thereby also giving an indication that there should exist a common proof
for both the finite and the infinitesimal case:
\begin{itemize} \item[(i)] we first use the (``$d\log$-like'') behaviour (cf. the comment
after (Definition \ref{dlog}))
$$\frac d{dc} \hatl_m\Big(c^\alpha (1-c)^\beta\Big) = \Big(\frac \alpha c -\frac \beta{1-c}\Big)
\hatl_{m-1}\Big(c^\alpha (1-c)^\beta\Big)\,$$
to group the terms of $\frac {d}{dc} \hatl_2\big( J(a,b,c)\big)$ with 
a coefficient $\frac 1 c$ (resp. $\frac 1 {1-c}$) together---these are exactly the terms with a
factor $c$ (resp. $1-c$). For instance, the terms with $\frac 1 c$ are as follows:
\begin{align*}
 \frac1 c \hatl_1\bigg(&(a-b+1)[ c] \\
&-a\Big[ \frac{c}{a}\Big] +
b\big[ \frac{c}{b}\Big] -b\Big[ \frac{ca}{b}\Big] -(1-b)\Big[ \frac{c(1-a)}{1-b}\Big] \\
&-c(1-a)\Big[ \frac{1-c^{-1}}{1-a^{-1}}\Big]  + c(1-b)\Big[ \frac{1-c^{-1}}{1-b^{-1}}\Big]  \\
&-c(a-b)\Big[ \frac{1-c^{-1}}{1-\frac{a}{b}}\Big]  -c(b-a)\Big[ \frac{1-c^{-1}}{1-\frac{1-a}{1-b}}\Big]  \bigg)
\end{align*}
In order to verify that this expression vanishes, we rewrite it in a slightly more convenient 
fashion (in order to be able to apply the four term relation line by line), neglecting the factor
$\frac1 c$, we get:
\begin{align*}\label{jred}
\,\hatl_1\Bigg(&(a-b+1)[ c] \\
&-a\left[ \frac{c}{a}\right] -c(1-a)\left[ \frac{1-c^{-1}}{1-a^{-1}}\right]  \\
&+b\left[ \frac{c}{b}\right] + c(1-b)\left[ \frac{1-c^{-1}}{1-b^{-1}}\right]  \\
&-b\left[ \frac{c}{\left(\frac{b}{a}\right)}\right] -c(a-b)\left[ \frac{1-c^{-1}}{1-\left(\frac{b}{a}\right)^{-1}}\right] \\
&-(1-b)\left[\frac{c}{\big(\frac{1-a}{1-b}\big)^{-1}}\right] -c(b-a)\left[ \frac{1-c^{-1}}{1-\frac{1-a}{1-b}}\right] \Bigg)\\
\end{align*}
Applying the 4-term equation (\ref{fund4}) ``linewise'' to the 2nd, 3rd, 4th and 5th line above 
with $x=a$, $x=b$, $x=\frac{b}{a}$ and $x=\frac{1-b}{1-a}$ , respectively,
this latter expression is  seen to reduce to
\begin{align*}
\hatl_1\Bigg(&(a-b+1)[ c] \\
& -a[ c] +c[ a] \ +b[ c] -c[ b] \\
&-a\left(\frac{b}{a}[ c] -
c\left[ \frac{b}{a}\right] \right)\ -\ (1-a)\left(\frac{1-b}{1-a}[ c] -
c\left[ \frac{1-b}{1-a}\right] \right)\Bigg)\\
=\hatl_1\Bigg(&c\left([ a] -[ b] +a\left[ \frac{b}{a}\right]  +(1-a)\left[ \frac{1-b}{1-a}\right] \right)\Bigg)
\end{align*}
which vanishes, again in view of the four term equation (and because the coefficients 
for $[c]$ add up to zero).
The terms with $\frac 1 {1-c}$ can be treated in a completely analogous way.\\
\item[(ii)]
The constant term in $c$ of the polynomial $\widehat{£_2}(J(a,b,c))$, i.e. the polynomial $\widehat{£_2}(J(a,b,0))$, is 
zero---this corresponds in the infinitesimal case to the degenerate case where we also put $c=0$ 
but where we need to give sense to expressions like $a[ \frac{b}{a}
]$ for $a=0$, the consistant choice being that it should be zero.
\\
\item[(iii)]
Instead of considering the coefficient of $c^p$ in the polynomial $\widehat{£_2}(J(a,b,c))$ we can equivalently
check that the constant coefficient in $c^p\widehat{£_2}(J(a,b,\frac1 c))$ is zero. In the infinitesimal case 
we can perform the same check using $c \widehat{£_2}(J(a,b,\frac1 c))$ (so we can use the analogy again).
\end{itemize}
\end{rem}



\bibliographystyle{plain}

\begin{thebibliography}{99}


\bibitem{ad} {\bf Acz\'el, J. and J. Dhombres} {\it Functional equations in several variables}, Encyclopedia of Math. and its Applications, Vol 31, Cambridge Univ. Press 1989.





\bibitem{Besser} {\bf Besser, A.} {\it Finite and $p$-adic polylogarithms}, preprint 2000.

\bibitem{bloch} {\bf Bloch, S.} {\it On the tangent space to Quillen $K$-theory}, Proc. Conf., Battelle
Memorial Inst., Seattle, Wash., 1972, Springer Lect. Notes in Math. n$ ^o$
341, (1973), pp. 205--210.   



\bibitem{bloch1} {\bf Bloch, S.} {\it Higher regulators, algebraic $K$-theory, and zeta-functions of elliptic curves}, Lecture Notes, U.C.~Irvine, 1977. 

\bibitem{bloch0} {\bf Bloch, S.} {\it
Applications of the dilogarithm function in algebraic $K$-theory and algebraic geometry}, Proc. of the
Int. Symp. on Alg. Geom. (Kyoto Univ., Kyoto, 1977), pp. 103--114, Kinokuniya Book Store, Tokyo, 1978. 







\bibitem{cathe-sl2} {\bf Cathelineau, J.-L.} {\it Sur l'homologie de SL$_{2}$ \`a coefficients dans l'action adjointe}, Math. Scand. 63 (1988), 51-86.


\bibitem{jlc} {\bf Cathelineau, J.-L.} {\em Homologie du groupe 
lin\'eaire et polylogarithmes (d'apr\`es Goncharov et d'autres)}, 
S\'eminaire Bourbaki, 772 (1992-93), Ast\'erisque 216 (1993) 311-341.

\bibitem{cathe-pol} {\bf Cathelineau, J.-L.} {\it Remarques sur les diff\'erentielles
des polylogarithmes uniformes}, Ann. Inst. Fourier, Grenoble, {\bf 46}, 5(1996), 1327-1347.

\bibitem{jlc-essen} {\bf Cathelineau, J.-L.} {\it Infinitesimal polylogarithms, multiplicative presentations of K\"ahler differentials and Goncharov complexes}, talk at the Workshop on Polylogarithms, Essen, May 1-4. (see \textsf{http://www.exp-math.uni-essen.de/$\sim$herbert/polyloquy.html})


\bibitem{coleman} {\bf Coleman, R.} 
{\it Dilogarithms, regulators and $p$-adic $L$-functions}, Invent. Math. 69
(1982), no. 2, 171--208.


\bibitem{dupont/sah} {\bf Dupont, J.~L.; Sah, C.-H.} 
{\it Scissors congruences. II.},
J. Pure Appl. Algebra 25,  (1982) 159-195.

\bibitem{pev} {\bf Elbaz-Vincent, Ph.} {\it The indecomposable K$_3$ of rings and homology of SL$_2$}, Jo. Pure. Appl. Alg. 132,(1998) 27-71.

\bibitem{pev-ktheory} {\bf Elbaz-Vincent, Ph.} {\it Homology of
Linear Groups with Coefficients in the Adjoint Action and $K$-theory},
$K$-Theory {\bf 16}, (1999), 35-50.

\bibitem{pev-essen} {\bf Elbaz-Vincent, Ph.} {\it Arithmetic and infinitesimal variations on a polylogarithm theme} (report on a work in progress with H. Gangl), talk at the Workshop on Polylogarithms, Essen, May 1-4. (see \textsf{http://www.exp-math.uni-essen.de/$\sim$herbert/polyloquy.html})

\bibitem{polyana2} {\bf Elbaz-Vincent, Ph.; Gangl, H.}  {\it On Poly(ana)logs II}, (in preparation).

\bibitem{aivpt} {\bf Elbaz-Vincent, Ph.; Gangl, H.}  {\it Arithmetical and Infinitesimal Variations on a Polylogarithmic Theme} (in preparation)


\bibitem{gangl1} {\bf Gangl, H.} {\it  Funktionalgleichungen von
    Polylogarithmen}, (German) Dissertation, Bonn, 1995. Bonner Mathematische
    Schriften {\bf 278}. Universität Bonn. 

\bibitem{gangl3} {\bf Gangl, H.} {\it Functional equations for higher logarithms},
    in preparation.


\bibitem{g1} {\bf Goncharov, A.B.} {\it Geometry of configurations, polylogarithms
and motivic cohomology}, Adv. Math. 114 (1995) 197-318.

\bibitem{g3} {\bf Goncharov, A.B.}  {\it Polylogarithms and motivic Galois groups},
Proc. of the Seattle conf. on motives, Seattle July 1991, 
AMS Proc. Symp. in Pure Math. 55 (1994) 2, 43-96.

\bibitem{gonmulpol} {\bf Goncharov, A.B.} {\it Multiple polylogarithms, cyclotomy and modular complexes}, Math. Res. Lett. 5 (1998), no. 4, 497--516. 





\bibitem{kontse} {\bf Kontsevich, M.} {\it The 1$\frac{1}{2}$-logarithm} (unpublished
note 1995).

\bibitem{lewin1} 
{\bf Lewin, L.} {\it Polylogarithms and associated functions} 
North-Holland Publishing Co., New York-Amsterdam, 1981. 

\bibitem{lewin2} {\bf Lewin, L. (Ed)} {\it Structural Properties of Polylogarithms}, Mathematical Surveys and Monographs, vol 37, AMS, 1991.

\bibitem{lewin3} {\bf Lewin, L.} {\it 
The order-independence of the polylogarithmic ladder structure---implications for a new category of functional
equations}, Aequationes Math. 30 (1986), no. 1, 1--20. 


\bibitem{oe} {\bf Oesterl\'e, J.} {\em Polylogarithmes}, S\'em. Bourbaki, 762
(1992-93), Ast\'erisque 216 (1993) 49-67. 




\bibitem{rib}  {\bf Ribenboim, P.} {\it 13 Lectures on Fermat's last theorem}, New York - Heidelberg - Berlin: Springer-Verlag.  (1979). 

\bibitem{Rog} {\bf Rogers, L.J.} {\it On function sum theorems connected with
    the series $\sum_1^\infty x^n/n^2$}, Proc London Math. Soc. {\bf 4}
    (1907), 169--189.

\bibitem{rosenberg} {\bf Rosenberg, J.} {\it Algebraic K-Theory and Its Applications}, GTM 147, Springer-Verlag, 1994. 



\bibitem{suslin} {\bf Suslin, A.A.} {\it $K\sb 3$ of a field, and the Bloch group. (Russian)} Galois theory, rings, algebraic groups and their applications (Russian). Trudy Mat. Inst. Steklov. 183
(1990), 180--199, 229.

\bibitem{washington} {\bf Washington, L.C.} {\it Introduction to cyclotomic fields}, GTM 83, Springer-Verlag 1997 (2nd edition).

\bibitem{weibel} {\bf Weibel, C.} {\it Private e-mail correspondence}, October 1999.

\bibitem{wojt1} {\bf Wojtkowiak, Z.} {\it A construction of analogs of the
    Bloch-Wigner function}, Math. Scand. {\bf 65}, 1989, p.140-142.

\bibitem{wojt2} {\bf Wojtkowiak, Z.} {\it A note on functional equations of
    the $p$-adic polylogarithms}, Bull. Soc. math. France {\bf 119}, 1991, p.343-370.

\bibitem{wojt3} {\bf Wojtkowiak, Z.} {\it The basic structure of polylogarithmic
functional equations}, Chapter 10 in  \cite{lewin2}, 205--231.


\bibitem{z1} {\bf Zagier, D.} {\it Polylogarithms, Dedekind zeta functions 
and the algebraic 
K-theory of fields}, Proc. Texel Conf. on Arithm. Alg. Geometry 1989, 
Birkh\"auser, 
Boston (1991) 391-430.

\bibitem{zag-app} {\bf Zagier, D.} {\it Special values and functional equations of polylogarithms},
Appendix A in \cite{lewin2}, 377--400.


\end{thebibliography}

\quad\\[5pt]
{\Small
\nobreak\indent{\sc Laboratoire G.T.A., UMR CNRS 5030, case courrier 51, Universit\'e Montpellier II,
  34095 Montpellier Cedex 5, France}.\\
\nobreak\indent{\itshape E-mail address} : \texttt{pev@math.univ-monpt2.fr}\\[3pt]
\nobreak\indent{\sc MPI Bonn, Vivatsgasse 7, D-53111 Bonn, Deutschland}.\\
\nobreak\indent{\itshape E-mail address} : \texttt{herbert@mpim-bonn.mpg.de}\\
\nobreak\indent{\itshape URL} :  \texttt{www.exp-math.uni-essen.de/\~{}herbert}
}
\newpage

\addcontentsline{toc}{part}{The $1{1\over 2}$-logarithm : Appendix by
Maxim Kontsevich}
\quad\\
\begin{center}
{\bf\Large  The $1{1\over 2}$-logarithm}\\[3pt] 
\sc\small (appendix to ``on poly(ana)logs i'' by ph.~elbaz-vincent and h.~gangl)\\[5pt] 
\normalsize maxim kontsevich\\[8pt]
\end{center}


{\sl  This appendix to the paper of Elbaz-Vincent and Gangl is included for historical purpose. It reproduces the text of \cite{kontse}, initially written for the private booklet ``Friedrich Hirzebruchs Emeritierung''.}\\[2pt]




Let $p> 2$ be a prime. Define a map from $\Z/p\Z$ to itself by the formula
$$H_p(x)=\sum_{k=1}^{p-1}  {x^k\over k}=x+{x^2\over 2}+\dots+{x^{p-1}\over p-1}\pmod p\,.$$
This function appears in explicit formulas for abelian extensions of cyclotomic fields. 
It looks like a truncated version of $\log(\frac{1}{1-x})$. Of course, it could not be a logarithm 
because there is no nonzero homomorphism from $(\Z/p\Z)^{\times}\simeq{\Z/(p-1)\Z}$ to $\Z/p\Z$. 
  I claim that $H_p$ is analogous to another well-known function of a real variable. I will derive the analogy 
by writing several functional equations for $H_p$. These equations will be  independent of $p$ and I will 
suppress the index $p$ from the notations.

(A):$\,\,\,\,\,H(1-x)=H(x)\,.$

\begin{proof} we can compute explicitly the coefficients of the polynomial $H(1-x)$.
  First of all, its zeroth coefficient is 
$H(1)=1+{1\over 2}+\dots{1\over p-1}=1+2+\dots +(p-1)={p(p-1)\over 2}=0\pmod p$.
  For $l$ between $1$ and $p-1$ the $l$-th coefficient of $H(1-x)$ is equal to
$$\sum_{k=1}^{p-1}\frac1k\,(-1)^l\,{k(k-1)\dots (k-l+1)\over l!}
={(-1)^l\over l!} \sum_{k=1}^{p-1}(k-1)\dots (k-l+1)=$$
$$=-\,\frac{(-1)^l}{l!}\,(0-1)(0-2)\dots(0-l+1)={(l-1)!\over l!}=\frac1l\;.$$
We use here the standard fact that
$$\sum_{k=0}^{p-1} P(k)=0$$ for any polynomial $P\in \Z/p\Z[x]$ of degree at most $p-2$.
\end{proof}
A simple generalization of the previous argument shows that\\

(B):$\,\,\,\,\, H(x+y)=H(y)+(1-y)\,H({x\over 1-y})+y\, H(-{x\over y})$ for $y\ne 0,\,1\,.$\\

Also there is a very elementary identity\\

(C):$\,\,\,\,\,\, x\,H({1\over x})= -H(x)$ for $x\ne 0$.\\

{\it Claim}: there is only one (up to a scalar factor) nonzero continuous solution of 
equations (A), (B), (C) in maps from $\R$ to itself. It is
 $$H_{\infty}(x)=-\left(x\,\,\log |x|+(1-x)\,\,\log |1-x|\right)\,.$$
Also the function $H_p$ is the unique (up to scalar factor) solution in maps from $\Z/p\Z$ to itself.\\[3pt]

\subsection*{Cohomological interpretation of functional equations}
Let $F$ be a field and suppose that $H:F\rightarrow F$ satisfies (A) and (B). The equation (C) will be irrelevant. 
We associate with $H$ a homogeneous function $\phi:F\times F\rightarrow F$ of degree~1:
$$\phi(x,y):=\begin{cases}(x+y) \,H\left({x\over x+y}\right)\,\qquad\text{if}\quad x+y\ne 0\,,\\
0\qquad\text{if}\quad x+y= 0\,.
\end{cases}$$

\noindent Equation (A) implies that $\phi(x,y)=\phi(y,x)$. Equation (B) is equivalent to the identity
$$\phi(x,y)-\phi(x,y+z)+\phi(x+y,z)-\phi(y,z)=0\,.$$
Thus, $\phi$ is a $2$-cocycle of the abelian group $F$ (the additive group of the field) with coefficients 
in itself as a trivial module. Because this cocycle is invariant under the usual action of the multiplicative 
 group $F^{\times}$ (acting both on the group and on the coefficients), we get a $2$-cocycle of the group 
of affine transformations of the line over $F$
$$\text{Aff}(1,F)=\{t\mapsto at+b| a\in K^{\times}, b\in K\}$$
with coefficients in the non-trivial $1$-dimensional representation given by the first coefficient. 
This $2$-cocycle defines an extension of $\text{Aff}(1,F)$ by $F$.
The resulting group $G$ can be identified as a set with $F\times F\times F^{\times}$.

Now we consider the case of $F=\R$ and assume that $H$ is a measurable map.
 There are no non-trivial measurable cohomology classes in $H^2(\R,\R)$, hence
 $\phi$ should be a coboundary. This means that there exists a function $\psi:\R
\rightarrow \R$ such that
$$\phi(x,y)=\psi(x)+\psi(y)-\psi(x+y)\,.$$
The homogeneity of $\phi$ implies that for any $\lambda\ne 0$ the function
$\psi_{\lambda}(x):=\psi(\lambda x)-\lambda \psi(x)$
 is additive in $x$. If we deal  with measurable maps only then $\psi_{\lambda}$ is a linear function. 
From this one can easily deduce that
$$\psi(x)/x=a\,\log |x| +b$$ 
for some $a,b\in \R$. Thus we get the solution  of functional equations for $F=\R$.

Now we turn to the case $F=\Z/p\Z$.  If the cohomology class in $H^2(F,F)\simeq \Z/p\Z$ corresponding 
to $H$ is zero then by arguments parallel to the previous one obtains a homomorphism from $(\Z/p\Z)^{\times}$ 
to $\Z/p\Z$. This homomorphism (a ``logarithm'') vanishes inevitably, thus giving the uniqueness of $H$ 
up to a scalar factor.

The group $G$ in the case of $F=\R$ is a $3$-dimensional solvable Lie group. The Lie algebra of $G$ 
is defined over $\Z$ and it has a base $x,y,z$ in which the commutation relations are
$$[x,y]=y,\,\,\,[x,z]=y+z,\,\,\,[y,z]=0\,.$$
This Lie algebra cannot be the Lie algebra of any algebraic group over $\Z$ or over $\Q$.
Nevertheless, we have defined groups of points  over $\R$ and over $\Z/p\Z$ for all odd primes~$p$.

\subsection*{Entropy}
The function $H$ is the entropy of a random variable taking two values. More generally, if $\xi$ takes a 
finite number of values with probabilities $p_1,\dots,p_k$, $\sum p_i=1$ then the entropy of $\xi$ is defined as
$$H(\xi):=-\sum_{i=1}^k p_i\,\,\log(p_i)\,.$$
We will consider the entropy also as a function of the collection of probabilities of elementary events, 
$H(\xi)=H(p_*)$.  
The main property of entropy is that if one random variable (say, $\xi$) is a function of another random 
variable (say, $\eta$) then the entropy of $\eta$ can be computed as follows. Let us denote probabilities 
of all possible values of $\eta$ by $p_{1,1},p_{1,2},\dots, p_{1,l_1};p_{2,1},\dots:\dots p_{k,l_k}$ 
in such a way that $p_{1,1}+p_{1,2}+\dots+p_{1,l_1}=p_1$ etc.
  Then we have $k$ conditional distributions of probabilities $p_{i,*}/p_i$ for each $i\le k$. The main 
identity of entropies is
$$H(p_{*,*})=H(p_*)+\sum_{i=1}^k p_i H({p_{i,*}\over p_i})\,,\,\,\,\,\,H(\eta)=H(\xi)+H_{\xi}(\eta)\,.$$
The last term in the formula above is the average value of the entropies of $\eta$ with given values 
of $\xi$ and it is called the relative entropy.

Using the main identity one can reduce by induction the calculation of the entropy of any random 
variable to the case of a two-valued variable, i.e. our function $H(x)$. One can check easily that the entropy 
of random variables computed using $H(x)$ is well-defined iff functional equations (A) and (B) are satisfied.

{\it Conclusion}:  If we have a random variable $\xi$ which takes finitely many values with all probabilities 
in $\Q$ then we can define not only the transcendental number $H(\xi)$ but also its ``residues modulo $p$"
for almost all primes $p\,\,$!

I propose calling the functions $H_p$ ``$1 {1\over 2}$-logarithms," because their functional equation 
contains $4$ terms, which is between 3 (the logarithm) and 5 (the dilogarithm giving an element in $H^3(Sl(2,\C),\R)$).

The natural question is to find functional equations for
the map $x\mapsto \sum_{k=1}^{p-1} x^k/k^2$ from $\Z/p\Z$ to itself. I don't know how to do it.

\quad\\[5pt]
{\Small
\nobreak\indent{\sc Maxim Kontsevich}\\
\nobreak\indent{\sc Institut des Hautes \'{E}tudes Scientifiques,  35, Route de
Chartres, 91440, Bures-sur-Yvette,  France}.\\
\nobreak\indent{\itshape E-mail address} : \texttt{maxim@ihes.fr}\\[3pt]
\end{document}